\documentclass{article}
\usepackage{amssymb,amsmath}
\usepackage{graphicx}

\newtheorem{theorem}{Theorem}[section]
\newtheorem{corollary}[theorem]{Corollary}
\newtheorem{example}[theorem]{Example}
\newtheorem{proposition}[theorem]{Proposition}
\newtheorem{lemma}[theorem]{Lemma}
\newtheorem{definition}[theorem]{Definition}
\newtheorem{remark}[theorem]{Remark}

\def\bD{\mathbb{D}}
\def\bR{\mathbb{R}}
\def\bZ{\mathbb{Z}}

\def\bV{\mathbb{V}}

\def\cB{\mathcal{B}}

\def\cD{\mathcal{D}}

\def\cF{\mathcal{F}}

\def\cL{\mathcal{L}}

\def\ve{\varepsilon}

\def\To{\longrightarrow}
\def\tostar{\stackrel{*}{\To}_S}
\def\indistes{\stackrel{*}{\To}_{\cD}}

\begin{document}

\title{Functional Convergence of Linear Sequences \\
in a non-Skorokhod Topology}

\author{Raluca Balan\footnote{Corresponding author. Department of Mathematics and Statistics, University of Ottawa,
585 King Edward Avenue, Ottawa, ON, K1N 6N5, Canada. E-mail address:
rbalan@uottawa.ca} \footnote{Research supported by a
grant from NSERC of Canada} \  Adam Jakubowski\footnote{Faculty of Mathematics and Computer Science, Nicolaus Copernicus University, 
87-100 Torun, Poland. E-mail address: adjakubo@mat.umk.pl
} \ and Sana Louhichi\footnote{Laboratoire Jean Kuntzmann, 
Institut de math\'ematiques appliqu\'ees de Grenoble, 51 rue des Mathématiques, 
F-38041 Grenoble, France. E-mail address: sana.louhichi@imag.fr  
}
}

\date{September 5, 2012}
\maketitle

\begin{abstract}
\noindent In this article, we prove a new functional limit theorem for the partial sum sequence $S_{[nt]}=\sum_{i=1}^{[nt]}X_i$ corresponding to a linear sequence of the form $X_i=\sum_{j \in \bZ}c_j \xi_{i-j}$ with i.i.d. innovations $(\xi_i)_{i \in \bZ}$ and real-valued coefficients $(c_j)_{j \in \bZ}$. This weak convergence result is obtained in space $\bD[0,1]$
endowed with the $S$-topology introduced in \cite{J97-EJP}, and the limit process is a linear fractional stable motion (LFSM). One of our result provides an extension of the results of \cite{avram-taqqu92} to the case when the coefficients $(c_j)_{j \in \bZ}$ may not have the same sign.
The proof of our result relies on the recent criteria for convergence in Skorokhod's $M_1$-topology (due to \cite{louhichi-rio11}), and a result which connects the weak $S$-convergence of the sum of two processes with the weak $M_1$-convergence of the two individual processes. Finally, we illustrate our results using some examples and computer simulations.
\end{abstract}

\noindent {\em MSC 2010 subject classification:} Primary 60F17; secondary 60G52

\section{Introduction}

The study of limit theorems for stochastic processes was initiated by Donsker in \cite{donsker51}   and Prohorov in \cite{prohorov54} in the case of processes with continuous trajectories, and continued by Skorokhod in
his seminal article \cite{skorokhod56}, in which he introduced the topologies $J_1, M_1, J_2, M_2$ on the space $\bD[0,1]$ of c\`adl\`ag functions on $[0,1]$. The basic idea is that criteria for compactness for sets in $\bD[0,1]$, once translated into criteria for tightness for probability measures on this space, become -via Prohorov's theorem- very powerful tools for proving functional convergence of stochastic processes, in the presence of finite-dimensional convergence. Although the immediate goal of Skorokhod's work was to extend the classical limit theory for sums of i.i.d. random variables to functional convergence, he saw this as part of a bigger program in which ``analysis of stochastic processes based on approximating them by simpler ones''  plays an important role. After the publication of Billingsley's cornerstone monograph \cite{billingsley68}, these ideas have been developed into a solid theory which
has been extended in many directions, and nowadays has ramifications in basically every area in probability theory.

In the present article, we study the functional convergence of the partial sum sequence $\{S_n=\sum_{i=1}^{n}X_i,n \geq 1\}$ corresponding to the linear sequence:
\begin{equation}
\label{linear-seq}
X_i=\sum_{j \in \bZ}c_j \xi_{i-j}, \quad i \in \bZ
\end{equation}
for suitably chosen coefficients $(c_j)_{j \in \bZ}$ and an i.i.d. sequence $(\xi_i)_{i \in \bZ}$. This problem has a very rich history and has been investigated by many authors. The philosophy behind these investigations is the one commonly encountered when dealing with a convolution between a random object (describing the shocks that drive the system), and a deterministic filter (the impulse-response function): if the filter is sufficiently smooth, then one expects that most properties of the input noise (in this case, the sequence $(\xi_i)_{i \in \bZ}$) can be transferred to the outcome result (in this case, the sequence $(X_i)_{i \in \bZ}$).

It is known that for a sequence $(\xi_i)_{i\in \bZ}$ of i.i.d. random variables, the class of distributional limits of its (suitably normalized) partial sum sequence coincides with the class of stable distributions, these being the only distributions which possess a domain of attraction. Here we used the common terminology according to which a random variable $\xi$ belongs to the domain of attraction of a distribution $G$ if there exist some constants $a_n>0$ and $b_n \in \bR$,
such that
\begin{equation}
\label{stable-LT}\frac{1}{a_n}\sum_{i=1}^{n}\xi_i-nb_n \To_{\cD} Z
\end{equation}
where $(\xi_i)_{i \geq 1}$ are i.i.d. copies of $\xi$ and $Z$ has distribution $G$.

The properties of random variables in the domain of attraction of stable distributions depend on the value of the index $\alpha \in (0,2]$ of stability, the case $\alpha=2$ corresponding to a normal limit distribution in (\ref{stable-LT}). Since the constants $a_n,b_n$ and the parameters of the distribution of $Z$ play an important role in the present article, we recall below their definitions. It is important to note that these objects are quite different in the case $\alpha=2$, versus $\alpha \in (0,2)$.

A comprehensive unified treatment which covers simultaneously the case $\alpha=2$ and $\alpha \in (0,2)$ is given by Feller in Section XVII.5 of \cite{feller71}. We describe very briefly the salient points. A (non-degenerate) random variable $\xi$ is in the domain of attraction of a normal distribution if and only if
$U(x)=E(\xi^2 1_{\{|\xi| \leq x\}})$ is slowly varying. In this case, the constants $(a_n)_{n \geq 1}$ are chosen such that
\begin{equation}
\label{def-an-alpha2}
nU(a_n) \sim c a_n^2
\end{equation}
for some constant $c>0$, $b_n=E(\xi)/a_n$ and $Z$ has a $N(0,c)$ distribution. Of course, when $E|\xi|^2<\infty$, $c={\rm Var}(\xi)$. On the other hand, a random variable $\xi$ is in the domain of attraction of a stable law of index $\alpha \in (0,2)$ if and only if
\begin{equation}
\label{RV-tail}
P(|\xi|>x) \sim x^{-\alpha}L(x) \quad \mbox{as} \ x \to \infty
\end{equation}
for some slowly varying function $L$, and
\begin{equation}
\label{balanced-tails}
\lim_{x \to \infty}\frac{P(\xi>x)}{P(|\xi|>x)}=p \quad \mbox{and} \quad
\lim_{x \to \infty}\frac{P(\xi<-x)}{P(|\xi|>x)}=q
\end{equation}
for some $p \geq 0$ and $q \geq 0$ with $p+q=1$. In this case, $b_n=0$ if $\alpha \in (0,1)$, $b_n=E(\xi)/a_n$ if $\alpha \in (1,2)$ and $b_n=E[\sin (\xi/a_n)]$ if $\alpha=1$.
The constants $(a_n)_{n}$ are chosen such that
\begin{equation}
\label{def-an-alpha}nP(|\xi|>a_n) \to C
\end{equation}
 for some constant $C>0$, and $Z$ has a stable $S_{\alpha}(\sigma,\beta,0)$ distribution, i.e.
$E(e^{iuZ})=e^{\Psi(u)}, u \in \bR$ with
$$\Psi(u)= \left\{
\begin{array}{ll}
-|u|^{\alpha} \sigma^{\alpha} \left(1-i {\rm sgn}(u) \beta \tan \frac{\pi \alpha}{2}\right)
& \mbox{if $\alpha \not=1$} \\
-|u| \sigma
\left(1+i {\rm sgn}(u) \beta \frac{2}{\pi} \ln |u|\right) & \mbox{if $\alpha=1$}
\end{array} \right.$$
and parameters $(\beta,\sigma)$ given by: $\beta=p-q$ and
\begin{equation}
\label{def-sigma}
\sigma^{\alpha}=\left\{
\begin{array}{ll}
C \frac{\Gamma(2-\alpha)}{1-\alpha}\cos(\frac{\pi \alpha}{2}) \ & \mbox{if $\alpha\not=1$} \\
C\frac{\pi}{2} & \mbox{if $\alpha=1$}
\end{array} \right.
\end{equation}

Following Skorokhod's ideas, the next step in this line of investigations
was to derive a functional analogue of (\ref{stable-LT}). In the case when $\xi$ is in the domain of attraction of a normal distribution, this is given by Donsker's theorem if  $E(\xi^2)<\infty$, respectively by  Proposition 1 of \cite{greenwood-resnick79} when $E(\xi^2)=\infty$. In general, such a result can be deduced from Theorem 2.7 of \cite{skorokhod57}. 
More precisely, using the same constants $a_n$ and $b_n$ as in (\ref{stable-LT}), one can prove that
\begin{equation}
\label{funct-stable}\frac{1}{a_n}\sum_{i=1}^{[n\cdot]}\xi_i-[n\cdot]b_n \To_{\cD} Z(\cdot)
\end{equation}
in the space $\bD[0,\infty)$ of c\`adl\`ag functions on $[0,\infty)$, equipped with the $J_1$-topology. The limit process $\{Z(t)\}_{t \geq 0}$ is a Brownian motion if $\alpha=2$, and an $\alpha$-stable L\'evy process if $\alpha \in (0,2)$. In both cases, $Z(1) \sim Z$.

Functional limit theorems for linear sequences of the form (\ref{linear-seq}) with $\sum_{j \in \bZ}|c_j|^2<\infty$ and (possibly dependent) innovations $(\xi_i)_{i \geq 1}$ with finite variance have received a lot of attention in the literature, usually by treating first the short-range dependence case (when $\sum_{i \geq 1}E(X_0 X_i)<\infty$), and then the more difficult case of long-range dependence (when $\sum_{i \geq 1}E(X_0 X_i)=\infty$). We refer the reader to \cite{davydov70}, \cite{peligrad-utev06}, \cite{DPP11} for a sample of relevant references. The remaining case when $\xi_1$ is in the domain of attraction of the normal law and has possibly infinite variance has been examined in the recent article \cite{peligrad-sang12}.

A totally new direction for the study of functional limit theorems for dependent sequences was initiated by Durrett and Resnick in \cite{durrett-resnick78}. In the i.i.d. case, this method supplied a surprising new proof of Skorokhod's result (\ref{funct-stable}) based on the convergence in distribution of the point process $N_n=\sum_{i \geq 1}\delta_{(i/n, \xi_i/a_n)}$ to the underlying Poisson process $N$ of the L\'evy process $\{Z(t)\}_{t \geq 0}$ (see Proposition 3.4 of \cite{resnick86}). The power of this method lies in the fact that it can be applied to very general dependent sequences, and supplies the convergence of a broad range of functionals (not only the sum). Since then, this method has been used by many authors in a variety of situations (e.g. \cite{davis-resnick85}, \cite{jakubowski-kobus89}, \cite{davis-hsing95}, \cite{davis-mikosch08}, \cite{tyran10}, \cite{BKS10}). In the case of linear sequences of the form (\ref{linear-seq}), a different method was used in \cite{phillips-solo92}, by expressing the sum $S_n=\sum_{i=1}^{n}X_i$ as the sum of an i.i.d. sequence and a negligible term.


From the results of \cite{astrauskas83} and \cite{davis-resnick85}, it can be inferred that:
\begin{equation}
\label{LT-linear}
\frac{1}{a_n}\sum_{i=1}^{[n\cdot]}X_i -A[n\cdot] b_n \stackrel{f.d.d.}{\To} A Z(\cdot), \quad \mbox{where} \quad A=\sum_{j \in \bZ}c_j,
\end{equation}
where $\stackrel{f.d.d.}{\To}$ denotes finite-dimensional convergence, $a_n,b_n$ and $Z(\cdot)$ are the same as in (\ref{funct-stable}), and the coefficients $(c_j)_{j \in \bZ}$ are such that $\sum_{j \in \bZ}|c_j| <\infty$ and
\begin{equation}
\label{cond-coeff}
\sum_{j \in \bZ}|c_j|^{\delta}<\infty \quad \mbox{for some} \quad 0<\delta<\alpha.
\end{equation}

Assume that $E(\xi)=0$ if $\alpha \in (1,2)$ and $\xi$ has a symmetric distribution if $\alpha=1$, so that $b_n=0$. A natural question is to see if (\ref{LT-linear}) can be extended to the functional convergence
\begin{equation}
\label{funct-LT-linear}\frac{1}{a_n}\sum_{i=1}^{[n\cdot]}X_i  \To_{\cD} A Z(\cdot)
\end{equation}
in the space $\bD[0,1]$ equipped with a suitable topology. In the well-known article \cite{avram-taqqu92}, Avram and Taqqu showed that if the series (\ref{linear-seq}) has at least two non-zero coefficients $c_j$, then (\ref{funct-LT-linear}) does not hold in the $J_1$-sense. Their argument showed that Skorokhod's criterion for $J_1$-tightness is not satisfied by a linear sequence (\ref{linear-seq}) with finitely many non-zero coefficients. However, in the case when all the coefficients are {\em non-negative}, they showed that (\ref{funct-LT-linear}) holds in the $M_1$-sense, assuming for $\alpha \in (1,2)$, that the coefficients satisfy a technical condition
(which holds for instance if $\delta<1$ and $(c_i)_{i \geq 0}$ and $(c_i)_{i<0}$ are monotone sequences).

The fact that (\ref{funct-LT-linear}) cannot hold in $\bD[0,1]$ equipped with {\em any} topology for which the supremum is continuous can be seen more easily by considering the linear sequence $$X_i=\xi_i-\xi_{i-1}, \quad i \in \bZ$$ whose coefficients are: $c_0=1$, $c_1=-1$ and $c_j=0$ for any $j \in \bZ \verb2\2 \{0,1\}$. In this case, the finite-dimensional distributions of $S_{[nt]}/a_n=(\xi_{[nt]}-\xi_0)/a_n,t \in [0,1]$ converge to $0$, but $\sup_{t \in [0,1]}S_{[nt]}/a_n=\max_{k \leq n}(\xi_n-\xi_0)/a_n$ converges to a random variable with a Fr\'echet distribution.

Finding a suitable topology on $\bD[0,1]$ for which (\ref{funct-LT-linear}) holds has remained an open problem in the literature since the article of Avram and Taqqu in 1992. In the present article, we give one possible answer to this question, using the $S$-topology introduced in \cite{J97-EJP} (which is weaker than $J_1$ and $M_1$).

Our first result (Theorem \ref{main-result}) is stated in the more general framework of \cite{KM88}, in which the normalizing constants for the partial sum sequence $S_n$ are of the form $d_n a_n$, for regularly varying constants $d_n$ of index $H-1/\alpha$ with $H \in [1/\alpha,1)$, and the limiting process is a {\em linear fractional stable motion} (LFSM), which is an $H$-self similar process defined as an integral with respect to $\{Z(t)\}_{t \in \bR}$.

Our second result (Theorem \ref{answer-open-problem}) is more clearly connected with the open problem mentioned above, and shows that if the coefficients $(c_j)_{j \in \bZ}$ admit a decomposition of the form
\begin{equation}
\label{decomp}
c_j=c_j'-c_j'', \quad \mbox{with} \quad c_j,c_j'' \geq 0
\end{equation}
for which both $(c_j')_{j \in \bZ}$ and $(c_j'')_{j \in \bZ}$ satisfy (\ref{cond-coeff})
for some $0<\delta<\alpha,\delta \leq 1$, then (\ref{funct-LT-linear}) holds in the space $\bD[0,1]$ equipped with the $S$ topology. By taking $c_j'=\max(c_j,0)$ and $c_j''=\max(-c_j,0)$, we see that this requirement is in fact equivalent to condition (\ref{cond-coeff}); for our purposes, it is more convenient to express (\ref{cond-coeff}) in this form.
In particular, (\ref{funct-LT-linear}) holds (in the $S$-sense) for a linear sequence (\ref{linear-seq}) with finitely many non-zero coefficients.

To prove this, we use the recent result of \cite{louhichi-rio11} which shows that for a strictly stationary sequence $(X_i)_{i \geq 1}$ which satisfies a dependence property called {\em association} (introduced in \cite{EPW67}), the finite dimensional convergence of the process $\{\gamma_n^{-1}\sum_{i=1}^{[nt]}X_i,t \in [0,1]\}$ to a c\`adl\`ag process $Y=\{Y(t)\}_{t \geq 0}$ is sufficient for its convergence in $\bD[0,1]$ equipped with $M_1$, provided that the normalizing constants $\gamma_n$ are regularly varying of index $\beta \in (0,2]$, and the tailsum of $Y(t)$ is of order $o(x^{-\beta/2})$. 
Using basic properties of the association, one observes that the linear sequence (\ref{linear-seq}) is associated, if either $c_j \geq 0$ for all $j \in \bZ$ or $c_j \leq 0$ for all $j \in \bZ$.

The second step in our argument is to add up the partial sums $S_{n}'$ and $S_{n}''$ corresponding to the linear sequences with coefficients $c_j'$ and $c_j''$, which arise from decomposition (\ref{decomp}). As shown in \cite{J97-EJP}, one of the remarkable properties of the $S$-topology is that addition is sequentially continuous. In the present article, we prove an analogue of this property for stochastic processes. More precisely, we show that for processes $(X_n)_n,(Y_n)_n,X,Y$ with trajectories in $\bD[0,1]$, if the finite dimensional distributions of $\{(X_n,Y_n)\}_n$ converge to those of $(X,Y)$ and both sequences $(X_n)_n$ and $(Y_n)_n$ converge in distribution in $\bD[0,1]$ equipped with the $M_1$-topology, then $(X_n+Y_n)_n$ converges in distribution to $X+Y$, in the space $\bD[0,1]$ equipped with the $S$ topology. This result provides us with the major ingredient needed to conclude our argument, since the final dimensional convergence of the (suitably normalize) pair $(S_n',S_n'')$ can be deduced without too much effort from the results of \cite{KM88}.

We conclude the introduction with few words about the organization of the present article. In Section \ref{section-S-top}, we recall the definition and main properties of the $S$-topology (as presented in \cite{J97-EJP}) and we prove the new result described above, regarding the convergence of a sum of two c\`adl\`ag processes. The proof of this result relies on the fact that the $S$-topology is weaker than $M_1$, a result whose proof we include here for the sake of completeness, since it was stated without proof in \cite{J97-EJP}. In Section \ref{section-result}, we recall the definition of the LFSM
and we derive our results regarding the convergence in $\bD[0,1]$ (equipped with the $S$-topology) of the partial sum corresponding to a linear sequence (\ref{linear-seq}) with coefficients $(c_j)_{j \in \bZ}$ with alternating signs. In Section \ref{section-examples}, we illustrate our results using some examples and simulations.
The Appendix contains some auxiliary results which are needed for the proof of the fact that the $S$-topology is weaker than $M_1$.

\noindent {\bf Acknowledgement:} The authors would like to thank Pierre-Yves Gaudreau-Lamarre for help with some of the simulations. His undergraduate project entitled ``Donsker theorem and its applications'' was funded by the Work-Study Program at the University of Ottawa, during the summer of 2012.

\section{The S-topology}
\label{section-S-top}

In this section, we recall the definition of the $S$-topology introduced in \cite{J97-EJP} and we prove some new properties, which are used in the present article.

Let $\bD[0,1]$ be the space of {\em c\`adl\`ag} functions $x:[0,1] \to \bR$, i.e. functions which are right-continuous on $[0,1)$ and have left-limits at each point $t \in (0,1]$. The space $\bD[0,1]$ can be endowed with the uniform topology, given by the norm: $$\|x\|_{\infty}=\sup_{t \in [0,1]}|x(t)|,$$ but also with the four Skorohod topologies ($J_1,M_1,J_2, M_2$) introduced in \cite{skorokhod56}.

The $S$-topology is a sequential topology on $\bD[0,1]$, defined using the concept of $\To_{S}$ convergence, which in turn relies on the concept of weak convergence for functions of bounded variation.

Let $\bV[0,1]$ be the space of functions $v:[0,1] \to \bR$ of bounded variation, endowed with the topology of weak-$*$ convergence: if $(v_n)_{n \geq 1}$ and $v$ are elements in $\bV[0,1]$, we write $v_n \To_{w} v$ if
$$\int_0^1 f(t) v_n(dt) \to \int_0^1 f(t)v(dt)$$
for any $f \in C[0,1]$, where $C[0,1]$ is the space of continuous functions on $[0,1]$.

Note that each function $v \in \bV[0,1]$ induces a signed measure $\nu_v$ on $[0,1]$ defined by $\nu_{v}([0,t])=v(t)$ for all $t \in [0,1]$. Since the space of the signed measures on $[0,1]$ can be identified with the dual space $(C[0,1])^*$ of $C[0,1]$ (endowed with the sup-norm metric), the space $\bV[0,1]$ can be endowed with the topology of weak-$*$ convergence of $(C[0,1])^*$.

The following result is crucial for the definition of the $\To_{S}$ convergence.

\begin{lemma}
\label{conv-in-V}
If $(v_n)_{n \geq 1}$ and $v$ are elements in $\bV[0,1]$ such that $v_n \To_{w} v$, then there exists a subsequence $(n_k)_{k \geq 1}$ and a countable set $D \subset [0,1)$ such that
$$v_{n_k}(t) \to v(t) \quad \mbox{for all} \ t \in [0,1] \verb2\2 D.$$
In addition, $v_n(1) \to v(1)$.
\end{lemma}

The proof of Lemma \ref{conv-in-V} follows by Banach-Steinhauss theorem and Helly's compactness theorem, using the fact that any $v \in \bV[0,1]$ can be written uniquely as $v=v(0)+v_1-v_2$ for some non-decreasing (c\`adl\`ag) functions $v_1,v_2$ with $v_1(0)=v_2(0)=0$. We omit the details.

\vspace{3mm}

We recall the definition of the $S$-convergence from \cite{J97-EJP}.

\begin{definition}
{\rm Let $(x_n)_{n \geq 1}$ and $x$ be elements in $\bD[0,1]$. We write $x_n \To_{S} x$ if for every $\ve>0$ there exist some functions $(v_{n,\ve})_{n \geq 1}$ and $v_{\ve}$ in $\bV$ such that:
\begin{eqnarray*}
& (i) & \|x_n -v_{n,\ve}\|_{\infty} \leq \ve \quad \mbox{for any} \ n \geq 1 \\
& (ii) & v_{n,\ve} \To_{w} v_{\ve}  \quad \mbox{for any} \ \ve>0 \\
& (iii) & \|x-v_{\ve}\|_{\infty} \leq \ve \quad \mbox{for any} \ \ve>0.
\end{eqnarray*}}
\end{definition}

The following result gives the explicit relationship between the sequence $(x_n)_{n \geq 1}$ and its limit $x$. Its proof is based on Lemma \ref{conv-in-V} and a diagonal argument. We omit the details.

\begin{lemma}
\label{def-S-limit}
Let $(x_n)_{n \geq 1}$ and $x$ be elements in $\bD[0,1]$. If $x_n \To_{S} x$, then
there exists a subsequence $(n_k)_{k \geq 1}$ and countable set $D \subset [0,1)$ such that
$$x_{n_k}(t)\to x(t) \quad \mbox{for all} \ t \in [0,1] \verb2\2 D.$$
In addition, $x_n(1) \to x(1)$.
\end{lemma}

The space $\bD[0,1]$ endowed with the $\To_S$ convergence is of type $\cL$, i.e.\\
(i) if $x_n=x,n \geq 1$ is a constant sequence, then $x_n \To_S x$; \\
(ii) if $x_n \To_S x$ then $x_{n_k} \To_S$ for any subsequence $(n_k)_k$.

Therefore, we can define the sequential topology induced by the $\To_S$ convergence. A set $F \subset \bD[0,1]$ is $S$-closed if for any $(x_n)_n \subset F$ with $x_n \To_S x$, we have $x \in F$. The $S$-topology is the collection of all sets $G \subset \bD[0,1]$ such that $G^c$ is $S$-closed.

We denote by $\tostar$ the convergence in the $S$-topology: $x_n \tostar x$ if for any $S$-open set $G$, there exists $N \geq 1$ such that $x_n \in G$ for all $n \geq N$. By the Kantorovich-Kisynski criterion, $x_n \tostar x$ if and only if for any subsequence $(n_k)_k$ there exists a further subsequence $(n_{k_l})_l$ such that $x_{n_{k_l}} \To_S x$.

\begin{remark}
\label{star-conv}
{\rm Note that, if $(x_n)_{n \geq 1}$ and $x$ are elements in $\bD[0,1]$ such that:\\
(i) $(x_n)_{n \geq 1}$ is $S$-relatively compact; and \\
(ii) for any subsequence $(n_k)_{k \geq 1}$ with $x_{n_k} \To_{S} y$ we have $y=x$,\\
then $x_{n} \tostar x$.}
\end{remark}

The most important property of the $S$-topology is the characterization of its relatively compact sets. This is expressed in terms of the number of upcrossings, or the number of oscillations, whose definitions we recall below.

For any real numbers $a<b$, let $N^{a,b}(x)$ be the {\em number of upcrossings}
of the interval $[a,b]$ by the function $x \in \bD[0,1]$, i.e. the largest integer $N \geq 1$ for which there exist some points
\begin{equation}
\label{points-def}
0 \leq t_1<t_2 \leq t_3<t_4 \leq \ldots \leq t_{2N-1}<t_{2N} \leq 1
\end{equation}
such that $x_{t_{2k-1}}<a$ and $x_{t_{2k}}>b$ for all $k=1,\ldots,N$.

For any $\eta>0$, let $N_{\eta}(x)$ be the {\em number of $\eta$-oscillations} of the function $x \in \bD[0,1]$, i.e. the largest integer $N \geq 1$ for which there exist some points $(t_i)_{1 \leq i \leq N}$ satisfying (\ref{points-def}), such that
$$|x(t_{2k})-x(t_{2k-1})|>\eta \quad \mbox{for all} \ k=1,\ldots,N.$$

Note that for any function $x \in \bD[0,1]$ and for any $a<b$,
$$N^{a,b}(x) \leq N_{b-a}(x)<\infty.$$

\begin{lemma}[Lemma 2.7 of \cite{J97-EJP}]
\label{charact-S-compact}
A set $A \subset \bD[0,1]$ is $S$-relatively compact if and only if it satisfies the following two conditions:
\begin{eqnarray*}
& (i) & \sup_{x \in A}\|x\|_{\infty}<\infty \\
& (ii) & \sup_{x \in A}N^{a,b}(x)<\infty \quad \mbox{for all} \ a,b \in \bR,a<b.
\end{eqnarray*}
Condition (ii) can be replaced by:\\
\begin{eqnarray*}
& (ii)' & \sup_{x \in A}N_{\eta}(x)<\infty \quad \mbox{for all} \ \eta>0.
\end{eqnarray*}
\end{lemma}

The relationship between the $S$-topology and the $M_1$-topology plays an important role in the present article. We refer the reader to the original Skorohod's article \cite{skorokhod56} for the definition of the $M_1$ topology, as well as Chapter 12 of \cite{whitt02} for a comprehensive account.

The $M_1$-convergence can be described using the {\em oscillation function}:
$$w(x,\delta):=\sup_{0 \vee (t_2-\delta) \leq t_1<t_2<t_3 \leq 1 \wedge (t_2+\delta)}H\big(x(t_1),x(t_2),x(t_3)\big),$$
for any $x \in \bD[0,1]$ and $\delta>0$, where $H(a,b,c)$ is the distance between $b$ and the interval with endpoints $a$ and $c$:
$$H(a,b,c)=(a\wedge c-a \wedge c \wedge b) \vee (a \vee c \vee b-a \vee c).$$

\begin{theorem}[2.4.1 of \cite{skorokhod56}]
\label{charact-M1-conv}
Let $(x_n)_{n \geq 1}$ and $x$ be arbitrary elements in $\bD[0,1]$. The following two conditions are necessary and sufficient for $x_n \to_{M_1} x$:
\begin{eqnarray*}
& (a) & x_n(t) \to x(t) \ \mbox{for any $t \in Q$, for a dense set $Q \subset [0,1]$ containing $0$ and $1$}; \\
& (b) & \lim_{\delta \to 0}\limsup_{n \to \infty}w(x_n,\delta)=0
\end{eqnarray*}
For the necessity part, $Q$ could be the set of continuity points of $x$, together with $0$ and $1$.
\end{theorem}

The following result was stated without proof in \cite{J97-EJP}. A short proof can be given using Skorohod's criterion 2.2.11 (page 267 of \cite{skorokhod56}) for the $M_1$-convergence, expressed in terms of the number of upcrossings. This proof has a clear disadvantage: it refers to an equivalent definition of the $M_1$-convergence,
which was not proved in Skorokhod's paper. In the present article, we give a new proof.

\begin{theorem}
\label{S-weaker-M1}
The $S$-topology is weaker than the $M_1$-topology (and hence, weaker than the $J_1$-topology).
Consequently, a set $A \subset \bD[0,1]$ which is $M_1$-relatively compact is also $S$-relatively compact.
\end{theorem}

\noindent {\bf Proof}: Let $(x_n)_{n \geq 1}$ and $x$ be elements in $\bD[0,1]$ such that $x_n \To_{M_1} x$. We will prove that
$x_n \tostar x$. For this, we apply Remark \ref{star-conv}. Suppose that $(n_k)_{k \geq 1}$ is an arbitrary subsequence such that $x_{n_k} \To_{S} y$. By Lemma \ref{def-S-limit}, there exists a subsequence $(n_{k_l})_l$ such that $x_{n_{k_l}}(t) \to y(t)$ for any $t \in Q_1$ for some dense set $Q_1 \subset [0,1]$ which contains $1$ (whose complement is countable). On the other hand, by Theorem \ref{charact-M1-conv}, $x_{n_k}(t) \to x(t)$ for any $t \in Q_2$ for some dense subset $Q_2 \subset [0,1]$ which contains $0$ and $1$ (whose complement is countable). Hence $y(t)=x(t)$ for any $t \in Q_1 \cap Q_2$, and $y=x$.

It remains to prove that $(x_n)_{n \geq 1}$ is $S$-relatively compact. By Lemma \ref{charact-S-compact}, it suffices to prove that:
\begin{eqnarray}
\label{uniform-bounded}
& & \sup_{n \geq 1}\|x_n\|_{\infty}<\infty \\
\label{oscilation-bounded}
& & \sup_{n \geq 1}N_{\eta}(x_n)<\infty \quad \mbox{for all} \ \eta>0.
\end{eqnarray}

Since $x_n \To_{M_1} x$ and the supremum is $M_1$-continuous, $\|x_{n}\|_{\infty} \to \|x\|_{\infty}$. Relation (\ref{uniform-bounded}) follows.

The remaining part of the proof is dedicated to (\ref{oscilation-bounded}). Let $\eta>0$ be arbitrary.

Let $0<\ve < \eta/2$. By Theorem \ref{charact-M1-conv}, there exist some $\delta > 0$ and an integer $n_0 \geq 1$ such that
\begin{equation}
\label{bound-osc-function}
w(x_n,\delta)
< \ve \quad \mbox{for all} \ n \geq n_0.
\end{equation}


Since the set $Q$ of continuity points of $x$ is dense in $[0,1]$, we can find some points $0 = t_0 < t_1 < \ldots <  t_{M} = 1$ in $Q$
such that for each $j=0, 1,\ldots,M-1$,
$$t_{j+1} - t_j < \delta.$$

By Theorem \ref{charact-M1-conv}, there exists an integer $n_1 \geq n_0$ such that for any $n \geq n_1$
\begin{equation}
\label{ineq-xn-x}| x_n(t_j) - x(t_j)| < \ve,\ \quad \mbox{for any} \ j = 0,1,\ldots, M.
\end{equation}

Fix an integer $n \geq n_1$. Suppose that there exist some points
$$ 0\leq s_1 < s_2 \leq s_3 < s_4 \leq \ldots \leq s_{2N-1} < s_{2N} \leq 1$$
such that
\begin{equation}\label{eaj3}
| x_n(s_{2k}) - x_n(s_{2k-1})| > \eta, \quad \mbox{for all} \ k=1,2, \ldots, N.
\end{equation}
The proof of (\ref{oscilation-bounded}) will be complete once we estimate the number $N$ by a constant independent of $n$.

Relation (\ref{eaj3}) says that the function $x_n$ has $N$ $\eta$-oscillations in the interval $[0,1]$. These oscillations can be divided into two (disjoint) groups. The first group (Group 1) contains the oscillations for which the corresponding interval $[s_{2k-1},s_{2k})$ contains at least one point $t_{j}'$. Since the number of points $t_j$ is $M$,
\begin{equation}
\label{group1}
\mbox{the number of oscillations in Group 1 is at most $M$}.
\end{equation}

In the second group (Group 2), we have those oscillations for which the corresponding interval $[s_{2k-1},s_{2k})$ contains no point $t_{j}$, i.e.
\begin{equation}\label{eaj4}
t_j \leq s_{2k-1} < s_{2k} \leq  t_{j+1} \quad \mbox{for some} \ j=0,1, \ldots,M-1.
\end{equation}



We now use Lemma \ref{lem-m1-number} (Appendix A).
Note that
$$\beta_n:=\sup_{t_{j} \leq u<v<w \leq t_{j+1}}H\big(x_n(u), x_n(v),x_n(w) \big) \leq w(x_n,\delta)<\ve,$$
and hence,
$$N_{\eta}(x_n,[t_j,t_{j+1}]) \leq \frac{2|x_n(t_{j+1})-x_n(t_{j})|+\beta_n}{\eta-\beta_n} <
\frac{ 4 K_{sup} + \varepsilon}{\eta-\varepsilon}$$
where $N_{\eta}(x_n,[t_j,t_{j+1}])$ is the number of $\eta$-oscillations of $x_n$ in the interval $[t_{j},t_{j+1}]$ and
$$K_{sup} =   \sup_{n \geq 1} \|x_n\|_{\infty} < +\infty.$$

Since there are $M$ intervals of the form $[t_j,t_{j+1}]$, we conclude that
\begin{equation}
\label{group2}
\mbox{the number of oscillations in Group 2 is at most}
\ M \cdot \frac{ 4 K_{sup} + \varepsilon}{\eta-\varepsilon}
\end{equation}

Using (\ref{group1}) and (\ref{group2}), we conclude that
$$N \leq   M\left(1+ \frac{ 4 K_{sup} + \varepsilon}{\eta-\varepsilon}\right)=M\frac{4K_{sup}+\eta}{\eta-\ve} ,$$
which does not depend on $n$. This concludes the proof of (\ref{oscilation-bounded}).
$\Box$

\vspace{3mm}

We now provide an example of a sequence $(x_n)_{n \geq 1}$ in $\bD[0,1]$ which is $S$-convergent, but does not converge in the $M_1$ topology.

\begin{example}
{\rm Let $x=0$ and
$$x_n(t)=1_{[1/2-1/n,1]}(t)-1_{[1/2+1/n,1]}(t)=\left\{
\begin{array}{ll} 1 & \mbox{if $\frac{1}{2}-\frac{1}{n} \leq t <\frac{1}{2}+\frac{1}{n}$} \\
0 & \mbox{otherwise}
\end{array} \right.$$
Then $x_n \To_{S} x$. To see this, we take $v_{n,\ve}=x_{n}$. Then $v_{n,\ve} \To_{w} v_{\ve}=0$ since for any $f \in C[0,1]$,
$$\int_{0}^{1}f(t)dv_n(t) =f\left(\frac{1}{2}+\frac{1}{n}\right)-f\left(\frac{1}{2}-\frac{1}{n}\right) \to 0.$$
The fact that $(x_n)_{n \geq 1}$ cannot converge in $M_1$ follows by Theorem \ref{charact-M1-conv} since if $t_1<\frac{1}{2}-\frac{1}{n}<t_2<\frac{1}{2}+\frac{1}{n}<t_3$, then $H\big(x_n(t_1),x_n(t_2), x_n(t_3)\big)=1$.}
\end{example}

We now consider random elements in $\bD[0,1]$ endowed with the $S$-topology. Recall that the Borel $\sigma$-field generated by the $S$-topology coincides with the Kolmogorov's $\sigma$-field $\cD$ generated by the projections $\pi_t,t \in [0,1]$
(see \cite{J97-EJP}). Therefore, a random element in $\bD[0,1]$ is a random variable $X:\Omega \to \bD[0,1]$ which is measurable with respect to $\cD$. Its law is a probability measure on $(\bD[0,1],\cD)$.

\vspace{2mm}

The following result is an immediate consequence of Theorem \ref{S-weaker-M1}.

\begin{corollary}
\label{tightness-S-M1}
Let $(P_{\alpha})_{\alpha \in \Lambda}$ be a family of probability measures on $\bD[0,1]$.
If $(P_{\alpha})_{\alpha \in \Lambda}$ is uniformly $M_1$-tight, then $(P_{\alpha})_{\alpha \in \Lambda}$ is also
uniformly $S$-tight.
\end{corollary}

By Lemma \ref{charact-S-compact}, one can give some criteria for the $S$-tightness of a family of probability measures on $\bD[0,1]$. More precisely, we have the following result.

\begin{proposition}[Proposition 3.1 of \cite{J97-EJP}]
\label{S-tight-criterion}
Let $(X_{\alpha})_{\alpha \in \Lambda}$ be a family of random elements in $\bD[0,1]$. The family $(X_{\alpha})_{\alpha \in \Lambda}$ is uniformly $S$-tight if and only if it satisfies the following two conditions:
\begin{eqnarray*}
& (i) & \{\|X_{\alpha}\|_{\infty}\}_{\alpha \in \Lambda} \ \mbox{is uniformly $S$-tight} \\
& (ii) & \{N^{a,b}(X_{\alpha})\}_{\alpha \in \Lambda} \ \mbox{is uniformly $S$-tight, for all $a,b \in \bR$, $a<b$}.
\end{eqnarray*}
Condition (ii) can be replaced by:
\begin{eqnarray*}
& (ii)' & \{N_{\eta}(X_{\alpha})\}_{\alpha \in \Lambda} \ \mbox{is uniformly $S$-tight, for all $\eta>0$}.
\end{eqnarray*}

\end{proposition}

Despite the fact that $\bD[0,1]$ equipped with the $S$-topology is not a metric space, the Direct Half of Prohorov's Theorem still holds, i.e. a uniformly $S$-tight family probability measures on $\bD[0,1]$ is relatively compact with respect to the $S$-weak convergence (defined using the convergence of integrals of $S$-continuous functions on $\bD[0,1]$). This follows by a strong form of Skorohod's Representation Theorem (Theorem 1.1 of \cite{J97-TPA}), using the fact that the $S$-topology possesses a countable family of $S$-continuous functions which separate the points in $\bD[0,1]$.

The Converse Part of Prohorov's Theorem also holds, but for this one needs to consider a stronger form of convergence in distribution (denoted by $\indistes$), given by the following definition (which was originally introduced in \cite{J00}).

\begin{definition}
\label{def-star-conv}
{\rm Let $(X_n)_{n \geq 1}$ be random elements in $\bD[0,1]$. We say that $(X_n)_{n \geq 1}$ is {\bf $\stackrel{*}{\To}_{\cD}$-convergent in distribution} if for any subsequence $(n_k)_k$ there exists a further subsequence $(n_{k_l})_l$ such that
\begin{center}
$(X_{n_{k_l}})_l$ admits a strong a.s. Skorohod representation,
\end{center}
i.e. there exist some random variables $(Y_l)_{l \geq 1}$ and $Y$ defined on $([0,1],\cB[0,1],\lambda)$ with values in $\bD[0,1]$ such that:\\
(i) $Y_l$ has the same distribution as $X_{n_{k_l}}$ for any $l \geq 1$;\\
(ii) $Y_l(\omega) \tostar Y(\omega)$ for any $\omega \in [0,1]$;\\
(iii) for any $\ve>0$ there exists an $S$-compact set $K \subset \bD[0,1]$ such that $$P(Y_l \in K \ \mbox{for all} \ l \geq 1)>1-\ve.$$}
\end{definition}

\begin{remark}
{\rm Let $g:\bD[0,1] \to \bR$ be an $S$-continuous bounded function. As a consequence of the previous definition, we obtain that $E[g(X_n)] \to E[g(Y)]$, since for any subsequence $(n_k)_k$ there exists a further subsequence $(n_{k_l})_l$ for which $E[g(X_{n_{k_l}})] \to E[g(Y)]$. This proves that $(\mu_n)_n$ converges weakly to $\mu$ (with respect to $S$), where $\mu_n$ is the law of $X_n$ and $\mu$ is the law of $Y$. If $X$ is a random element in $\bD[0,1]$ with law $\mu$, we write $X_n \indistes X$. }
\end{remark}

The following two results are needed for the proof of Theorem \ref{sum-conv} below. We recall them for the sake of completeness.

\begin{theorem}[Theorem 3.4 of \cite{J97-EJP}]
\label{theorem3-4}
Let $(X_{\alpha})_{\alpha \in \Lambda}$ be a family of random elements in $\bD[0,1]$. Then $(X_{\alpha})_{\alpha \in \Lambda}$ is uniformly $S$-tight if and only if it is relatively compact with respect to $\indistes$.
\end{theorem}

\begin{theorem}[Theorem 3.5 of \cite{J97-EJP}]
\label{theorem3-5}
Let $(X_n)_{n \geq 1}$ and $X$ be random elements in $\bD[0,1]$ such that:
\begin{eqnarray*}
&(i) & (X_n(t_1), \ldots,X_n(t_k)) \To_{\cD} (X(t_1),\ldots,X(t_k)) \ \mbox{for any $t_1<\ldots<t_k$ in $Q$}, \\
& & \mbox{and $k \geq 1$, where $Q$ is a dense set in $[0,1]$ which contains $1$}; \\
&(ii) & (X_n)_{n \geq 1} \ \mbox{is relatively compact with respect to $\indistes$}.
\end{eqnarray*}
Then $X_n \indistes X$.
\end{theorem}

We conclude this section with a new result which shows that the sum of two processes converges in distribution (in the sense of $\indistes$) provided that the two processes converge in $M_1$ and the finite-dimensional distributions converge. This result will be used in Section \ref{section-result}.

\begin{theorem}
\label{sum-conv}
Let $(X_n)_{n \geq 1}$, $(Y_n)_{n \geq 1}$, $X$ and $Y$ be random elements in $\bD[0,1]$ such that:
\begin{eqnarray*}
& (i) & X_n \To_{\cD} X \ \mbox{in $D[0,1]$ equipped with $M_1$}; \\
& (ii)& Y_n \To_{\cD} Y \ \mbox{in $D[0,1]$ equipped with $M_1$};\\
& (iii) & (X_n(t_1),Y_n(t_1), \ldots, X_n(t_k),Y_n(t_k)) \To_{\cD} (X(t_1),Y(t_1), \ldots, X(t_k),Y(t_k)) \\
 & & \mbox{for any $0 \leq t_1<\ldots<t_k \leq 1$ and $k \geq 1$}.
\end{eqnarray*}
Then $X_n+Y_n \indistes X+Y$.
\end{theorem}

\noindent {\bf Proof:} First, let us observe that both sequences $(X_n)_{n \geq 1}$ and $(Y_{n})_{n \geq 1}$ are uniformly $M_1$-tight. This follows by Le Cam's theorem (Theorem 8 in Appendix III of \cite{billingsley68}), since a single probability measure on $\bD[0,1]$ is $M_1$-tight.

By Corollary \ref{tightness-S-M1}, both sequences $(X_n)_{n \geq 1}$ and $(Y_{n})_{n \geq 1}$ are uniformly $S$-tight, and hence, they satisfy conditions $(i)$ and $(ii)'$ of Proposition \ref{S-tight-criterion}. Since for any functions $x,y \in \bD[0,1]$,
$\|x+y\|_{\infty} \leq \|x\|_{\infty}+\|y\|_{\infty}$
and
$$N_{\eta}(x+y)\leq N_{\eta/2}(x)+N_{\eta/2}(y) \quad \mbox{for any} \ \eta>0,$$
it follows that the sequence $(X_n+Y_n)_{n \geq 1}$ also satisfies conditions $(i)$ and $(ii)'$ of Proposition \ref{S-tight-criterion}. Hence, the sequence $(X_n+Y_n)_{n \geq 1}$ is uniformly $S$-tight. By Theorem \ref{theorem3-4}, $(X_n+Y_n)_{n \geq 1}$ is relatively compact with respect to $\indistes$.

Finally, by our hypothesis $(iii)$ and the Continuous Mapping Theorem,
$$\big(X_n(t_1)+Y_n(t_1),\ldots,X_n(t_k)+Y_n(t_k)\big) \To_{\cD} \big(X(t_1)+Y(t_1),\ldots, X(t_k)+Y(t_k)\big)$$
for any $0 \leq t_1<\ldots<t_k \leq 1$ and $k \geq 1$.
The result follows by applying Theorem \ref{theorem3-5} to the sequence $(X_n+Y_n)_{n \geq 1}$. $\Box$

\section{Our Results}
\label{section-result}

In this section, we derive some new results regarding the convergence in distribution (in the $S$-topology) of the partial sum of the linear sequence (\ref{linear-seq}).

Let $(X_i)_{i \in \bZ}$ be the linear sequence given by (\ref{linear-seq}), where $(\xi_i)_{i \in \bZ}$ are i.i.d. random variables in the domain of attraction of a stable law of index $\alpha \in (0,2]$, and the coefficients $(c_j)_{j \in \bZ}$ satisfy (\ref{cond-coeff}).

As in \cite{avram-taqqu92}, we assume that $E(\xi_1)=0$ if $\alpha \in (1,2]$ and $\xi_1$ has a symmetric distribution if $\alpha=1$. Note that the series appearing in the right-hand side of (\ref{linear-seq}) converges a.s. This follows by applying Theorem 12.10.4 of \cite{kawata72} (when (\ref{cond-coeff}) holds with $\delta>1$), or Theorem 12.11.2 of \cite{kawata72} (when (\ref{cond-coeff}) holds with $\delta \leq 1$).

For any $n \geq 1$, we consider the partial sum processes:
$$Z_n(t)=\frac{1}{a_n}\sum_{i=1}^{[nt]}\xi_i,\ t \geq 0 \quad \mbox{and} \quad S_n(t)=\sum_{i=1}^{[nt]}X_i,\ t\geq 0.$$

Due to our assumptions, relation (\ref{funct-stable}) holds with $b_n=0$, i.e.
$$Z_n(\cdot) \To_{\cD} Z(\cdot)$$ in $\bD[0,\infty)$ equipped with the $J_1$-topology.
The constants $(a_n)_n$ and the distribution of $Z(1) \sim Z$
have been specified in the introduction.

\begin{remark}
\label{def-Z-process}
{\rm For any $\alpha \in (0,2]$, $Z$ has a {\em strictly} $\alpha$-stable distribution.
When $\alpha=2$, $\{Z(t)\}_{t \geq 0}$ is a Brownian motion with variance $c$. In this case, $Z(t)-Z(s)$ has a $N(0,c(t-s))$ distribution, for any $s<t$.

When $\alpha \in (0,2)$, $\{Z(t)\}_{t \geq 0}$ is an $\alpha$-stable L\'evy motion, in the sense of Definition 3.1.3 of \cite{ST94} (except that in this definition, it is required that $\sigma=1$). To see this, note that for any $s<t$, by property  1.2.3 (page 11 of \cite{ST94}) 
$$Z(t)-Z(s) \sim Z(t-s) \sim (t-s)^{1/\alpha} Z \sim
\left\{
\begin{array}{ll}
S_{\alpha}(\sigma(t-s)^{1/\alpha},\beta,0) & \mbox{if $\alpha \not=1$} \\
S_{\alpha}(\sigma(t-s)^{1/\alpha},\beta,\mu_1) & \mbox{if $\alpha =1$}
\end{array} \right.$$
where $\mu_1=-\frac{2\sigma \beta}{\alpha \pi} (t-s)^{1/\alpha} \ln (t-s)$.
If $\alpha=1$, we assume that $p=q$ and $\mu_1=0$.


}
\end{remark}

We follow very closely the approach of \cite{KM88}. We first observe that $a_n^{-1}S_n(t)$ can be expressed as an integral with respect to the process $Z_n=\{Z_n(t)\}_{t \geq 0}$.
Recall that the integral of a deterministic function $f:\bR \to \bR$ with respect to $Z_n$ is defined by:
\begin{equation}
\label{def-int-Zn}
\int_{-\infty}^{\infty}f(u)dZ_n(u):=\frac{1}{a_n}\sum_{j \in \bZ}f(j/n)\xi_j,
\end{equation}
provided that the sum converges a.s. One sufficient condition for this is:
\begin{equation}
\label{conv-sum-f}\sum_{j \in \bZ} |f(j/n)|^{b}<\infty \quad \mbox{for some} \quad 0<b<\alpha.
\end{equation}

In our case, using (\ref{linear-seq}) and Fubini's theorem, we obtain that
$$S_n(t)=\sum_{i=1}^{[nt]}\sum_{j \in \bZ}c_{i-j}\xi_j=\sum_{j \in \bZ} \left(\sum_{k=1-j}^{[nt]-j}c_k\right)\xi_j,$$
and hence
$$\frac{1}{a_n}S_n(t)=\frac{1}{a_n}\sum_{j \in \bZ}f_n\left(t,\frac{j}{n}\right)\xi_j=\int_{-\infty}^{\infty}f_n(t,u)dZ_n(u),$$ where $$f_n(t,u)=\sum_{k=1-[nu]}^{[nt]-[nu]} c_k, \ u \in \bR.$$

\begin{remark}
{\rm Note that condition (\ref{conv-sum-f}) satisfied for $f=f_n(t,\cdot)$. This follows by taking $b=\delta$ since
\begin{eqnarray*}
\sum_{j \in \bZ} \left|\sum_{k=1-j}^{[nt]-j}c_k\right|^{\delta} &\leq & C_{n,t,\delta} \sum_{j \in \bZ}\sum_{k=1-j}^{[nt]-j}|c_k|^{\delta} =C_{n,t,\delta} \sum_{k \in \bZ} \sum_{j=1-k}^{[nt]-k}|c_k|^{\delta}\\
&=& C_{n,t,\delta} [nt] \sum_{k \in \bZ}|c_k|^{\delta}<\infty
\end{eqnarray*}
where $C_{n,t,\delta}=1$ if $\delta\leq 1$ and $C_{n,t,\delta}=[nt]^{\delta-1}$ if $\delta>1$.

}
\end{remark}

In \cite{KM88}, it is shown that the finite dimensional distributions of the process $S_n/(d_n a_n)$ converge to those of a LFSM, for a sequence $(d_n)_n$ of suitable constants. We will need this result below. To introduce the LFSM in the case $\alpha \in (0,2)$, we have to recall the definition of an $\alpha$-stable random measure.

\begin{definition} {\rm Let $m$ be a positive measure on $(\bR,\cB)$ and $\cB_0=\{A \in \cB;m(A)<\infty\}$. Let $\alpha \in (0,2)$. A collection $M=\{M(A);A \in \cB_0\}$ of random variables defined on a probability space $(\Omega,\cF,P)$ is called an {\em $\alpha$-stable random measure} on $(\bR,\cB)$ with control measure $m$ and skewness intensity $\beta \in [-1,1]$ if: \\
(i) 
for any disjoint sets $A_1, \ldots, A_k \in \cB_0$, $M(A_1), \ldots, M(A_k)$ are independent;\\
(ii) 
for any disjoint sets $(A_n)_{n \geq 1} \subset \cB_0$ with $\cup_{n \geq 1}A_n \in \cB_0$, $M(\cup_{n \geq 1}A_n)=\sum_{n \geq 1}M(A_n)$ a.s.;\\
(iii) for any $A \in \cB_0$, $M(A) \sim S_{\alpha}(m(A)^{1/\alpha}, \beta,0)$. }
\end{definition}

The existence of $M$ is shown in Section 3.3 of \cite{ST94}. The stochastic integral
$$I(f)=\int_{\bR}f(u)M(du)$$
is defined in Section 3.4 of \cite{ST94} for any measurable function $f:\bR \to \bR$ such that $\int_{\bR}|f(u)|^{\alpha}m(du)<\infty$ if $\alpha \not=1$, and $\int_{\bR} |f(u) \ln |f(u)||m(du)<\infty$ if $\alpha=1$.
By Property 3.2.2 (page 117 of \cite{ST94}),
$$I(f) \sim S_{\alpha}(\sigma_f,\beta_f,\mu_f)$$
where $\sigma_f^{\alpha}=\int_{\bR} |f(u)|^{\alpha}m(du)$, $\mu_f=-\frac{2\beta}{\pi} \int_{\bR} f(u) \ln |f(u)| m(du)$ if $\alpha=1$,
\begin{equation}
\label{def-mu-f}
\mu_f=0 \quad \mbox{if}  \quad \alpha \not=1 \quad \mbox{and} \quad \beta_f=\frac{\beta\int |f(u)|^{\alpha} {\rm sgn}f(u)m(du)}{\int |f(u)|^{\alpha} m(du)}.
 \end{equation}

In what follows, we assume that $m=\sigma^{\alpha}\lambda$ where $\lambda$ is the Lebesgue measure on $\bR$ and $\sigma$ is given by (\ref{def-sigma}). Then
$$M((s,t]) \sim S_{\alpha}(\sigma (t-s)^{1/\alpha},\beta,0) \quad \mbox{for any} \quad s<t.$$

Setting $Z_M(t):=M([0,t]),t \in \bR$, we see that $\{Z_M(t)\}_{t \geq 0}$ has the same finite-dimensional distributions as $\{Z(t)\}_{t \geq 0}$. For this reason, we say that $Z_M$ is an extension of $Z$ to $\bR$. We denote $Z_M$ simply by $Z$ and $\int_{\bR} f(u)M(du)$ by $\int_{\bR} f(u)dZ(u)$. We are now ready to give the definition of the LFSM.

\begin{definition}[Definition 7.4.1 of \cite{ST94}]
{\rm The {\em linear fractional stable motion} (LFSM) is the stochastic process
$\{\Lambda_{\alpha,H,a,b}(t)\}_{t \in \bR}$ given by:
$$\Lambda_{\alpha,H,a,b}(t)=\int_{-\infty}^{\infty}f_{\alpha,H,a,b}(t,u)dZ(u),$$
where $\alpha \in (0,2]$, $H \in (0,1)$, $H \not=1/\alpha$, $a\in \bR$, $b \in \bR$ with $|a|+|b|>0$ and
$$f_{\alpha,H,a,b}(t,u)=a\{(t-u)_+^{H-1/\alpha}-(-u)_+^{H-1/\alpha}\}+
b\{(u-t)_+^{H-1/\alpha}-(u)_+^{H-1/\alpha}\}.$$}
\end{definition}

\begin{remark}
\label{path-properties}
{\rm (i) The {\em fractional Brownian motion} (FBM) of Hurst index $H \in (0,1)$ is a zero-mean Gaussian process $\{B_{H}(t)\}_{t \geq 0}$ with  $E[B_H(t)B_H(s)]=(t^{2H}+s^{2H}-|t-s|^{2H})/2$. This process admits the ``moving average'' representation: (see Proposition 7.2.6 of \cite{ST94})
$$B_H(t)=C_H \int_{-\infty}^{\infty}\{(t-u)_+^{H-1/2}-(-u)_+^{H-1/2}\}dB(u)$$
where $\{B(t)\}_{t \in \bR}$ is a two-sided standard Brownian motion and
\begin{equation}
\label{definition-CH}
C_{H}=\left\{ \int_{0}^{\infty} \left[(1+u)^{H-1/2}-u^{H-1/2}\right]^2du+\frac{1}{2H} \right\}^{1/2}.
\end{equation}
Therefore, the process $\Lambda_{2,H,c_H,0}$ is a FBM. Note that the sample paths of the FBM are $\gamma$-H\"older continuous, for any $\gamma \in (0,H)$.

(ii) It is convenient to have a unified notation which covers also the case $H= 1/\alpha$. Therefore, if $H=1/\alpha$ we let $f_{\alpha,H,a,b}(t,u)=(a-b)1_{(0,t]}(u)$ and
$$\Lambda_{\alpha,H,a,b}(t)=\int_{-\infty}^{\infty}f_{\alpha,H,a,b}(t,u)dZ(u)=(a-b)Z(t).$$
In this case, $\Lambda_{\alpha,H,a,b}$ has a c\`adl\`ag modification (see e.g. Theorem 5.4 of \cite{resnick07}).

(iii) Suppose that $\alpha \in (1,2)$ and $H>1/\alpha$. By Proposition 7.4.2 of \cite{ST94}, $\Lambda_{\alpha,H,a,b}$ is an $H$-sssi process, i.e. it is $H$-self similar and has stationary increments.
Note that $$\Lambda_{\alpha,H,a,b}(t)=I(f) \sim S_{\alpha}(\sigma_f,\beta_f,\mu_f)$$ with $f=f_{\alpha,H,a,b}$ and parameters $(\sigma_f,\beta_f,\mu_f)$ as above. By (\ref{def-mu-f}) and our assumption that $\beta=0$ if $\alpha=1$, it follows that $\Lambda_{\alpha,H,a,b}(t)$ has a strictly $\alpha$-stable distribution.
By Theorem 12.4.1 of \cite{ST94}, $\Lambda_{\alpha,H,a,b}$ has a continuous version.
}
\end{remark}

The following recent result of \cite{louhichi-rio11} lies at the origin of our investigations. We recall it for the sake of completeness.

\begin{theorem}[Theorem 1 of \cite{louhichi-rio11}]
\label{LR-theorem}
Let $(X_i)_{i \geq 1}$ be a strictly stationary sequence of associated random variables and $S_n(t)=\sum_{i=1}^{[nt]}X_i, t \in [0,1]$. Let $(\gamma_n)_{n \geq 1}$ be a non-decreasing sequence of constants such that $\gamma_n \to \infty$ and
$$\gamma_n \sim n^{1/\beta}L(n)$$
for some $\beta \in (0,2]$ and a slowly varying function $L$.
Let $(b_n)_{n \geq 1}$ be a sequence of real numbers such that $K:=\sup\{|b_k-b_n|;n \geq 1, n \leq k \leq 2n \}<\infty$. Assume that either $K=0$ or $\liminf_{n}\gamma_n/n>0$. Let
$$Y_n(t)=\frac{1}{\gamma_n}(S_n(t)-[nt]b_n) \quad t \in [0,1].$$ If $Y_n(\cdot) \stackrel{f.d.d.}{\To} Y(\cdot)$ where $Y=\{Y(t)\}_{t \in [0,1]}$ is a c\`adl\`ag process with
\begin{equation}
\label{tailsum-cond}
\lim_{x \to \infty}x^{\beta/2}P(|Y(t)| \geq x)=0 \quad \mbox{for all} \quad t \in [0,1],
\end{equation}
then $Y_n(\cdot) \To_{\cD} Y(\cdot)$ in $(\bD[0,1],M_1)$.

\end{theorem}

Based on Theorem \ref{LR-theorem}, we can derive our first result, which is an extension of Theorem 5.1 of \cite{KM88} from the convergence of finite-dimensional distributions to the convergence in distribution in $\bD[0,1]$, in the $M_1$-sense (in the case of non-negative coefficients). For this, we need to observe that
\begin{equation}
\label{rate-an}
a_n \sim n^{1/\alpha}L_1(n)
\end{equation}
for a slowly varying function $L_1$.

\begin{theorem}
\label{th-pos-coeff}
Let $\alpha \in (1,2]$. Suppose that $c_j \geq 0$ for all $j \in \bZ$ and
$$\frac{1}{d_n}\sum_{j=0}^{n}c_j \to a, \quad \frac{1}{d_n}\sum_{j=-n}^{0}c_j \to -b, \quad  c_n=O\left(\frac{d_{|n|}}{|n|}\right),$$
for some $a,b \in \bR$ with $|a|+|b|>0$, where
\begin{equation}
\label{rate-dn}
d_n=n^{H-1/\alpha}L_2(n)
\end{equation} for some $1/\alpha \leq H<1$ and a slowly varying function $L_2$. Let $(a_n)_n$ be chosen to satisfy (\ref{def-an-alpha2}) or (\ref{def-an-alpha}) depending on whether $\alpha=2$ or $\alpha \in (0,2)$. Then
$$\frac{1}{d_n a_n} S_n(\cdot) \To_{\cD} \Lambda_{\alpha,H,a,b}(\cdot) \quad \mbox{in} \quad (\bD[0,1],M_1).$$
\end{theorem}

\noindent {\bf Proof:} To simplify the notation we write $\Lambda$ instead of $\Lambda_{\alpha,H,a,b}$. By Remark \ref{path-properties}, we may assume that $\Lambda$ is a c\`adl\`ag process.
We apply Theorem \ref{LR-theorem} with $b_n=0$, $\gamma_n=d_n a_n$ and $\beta=1/H$. By Theorem 5.1 of \cite{KM88},
$$\frac{1}{d_n a_n}S_n(\cdot) \stackrel{f.d.d.}{\To} \Lambda(\cdot).$$

\noindent Due to (\ref{rate-an}) and (\ref{rate-dn}), $c_n \sim n^{H}L(n)$ for the slowly varying function $L:=L_1 L_2$. If $\alpha=2$, $\Lambda(t)$ has a normal distribution,
whereas if $\alpha \in (0,2)$, $\Lambda(t)$ has an $\alpha$-stable distribution.
In both cases, (\ref{tailsum-cond}) holds.

Finally, we prove that $(X_i)_{i \in \bZ}$ is an associated sequence. By definition, $X_i=\lim_{N \to \infty}X_i^{(N)}$ a.s. where $X_i^{(N)}=\sum_{j=-N}^{N}c_{j}\xi_{i-j}$. By property (P5) of \cite{EPW67}, it suffices to show that $(X_i^{(N)})_{i \in \bZ}$ is associated for any $N \geq 1$. Note that $X_i^{(N)}=f(\xi_{i-N}, \ldots, \xi_{i+N})$, where the function $f: \bR^{2N+1} \to \bR$ is defined by $f(x_{i-N}, \ldots,x_{i+N})=\sum_{j=-N}^{N}c_j x_{i-j}$. Since $c_j \geq 0$ for all $j \in \bZ$, the function $f$ is coordinate-wise non-decreasing.
Let $i_1< \ldots<i_k$ be a finite set of indices in $\bZ$. Let $I=\{i_1-N, \ldots, i_k+N\}$. Then for each $j \in \{1, \ldots,k\}$, we can say that $X_{i_{j}}^{(N)}=f_j(\xi_i; i \in I)$ for some coordinate-wise non-decreasing function $f_j$. Since $(\xi_{i})_{i \in I}$ are associated (Theorem 2.1 of \cite{EPW67}), by property (P4) of \cite{EPW67}, it follows that $\{X_{i_1}^{(N)}, \ldots, X_{i_k}^{(N)}\}$ are associated. $\Box$

\begin{remark}
{\rm (i) Under the assumptions of Theorem \ref{th-pos-coeff}, if $H>1/\alpha$
then
$$\frac{1}{d_n a_n}S_n(\cdot) \To_{\cD} \Lambda_{\alpha,H,a,b}(\cdot) \quad \mbox{in} \quad (\bD[0,1],U)$$
since $\{\Lambda_{\alpha,H,a,b}(t)\}_{t \in [0,1]}$ is continuous. Here $U$ denotes the uniform topology. When $H=1/\alpha$, Theorem \ref{th-pos-coeff} can be viewed as a variant of Theorem 2 of \cite{avram-taqqu92}; in this case, the limit is $(a-b)Z(\cdot)$.

(ii) In the case $H<1/\alpha$, Theorem 5.2 of \cite{KM88} gives the finite-dimensional convergence of $(a_n^{-1}S_n(\cdot)-A Z_n(\cdot))/d_n$ where $A=\sum_{j \in \bZ}c_j$. However, this process cannot be identified with the partial sum process corresponding to a sequence of associated random variables. Therefore, the case $H<1/\alpha$ cannot be treated by the methods of the present article. }
\end{remark}


The following theorem is the main result of this article.

\begin{theorem}
\label{main-result}
Let $\alpha \in (1,2]$. Suppose that $c_j=c_j'-c_j''$ for some $c_j',c_j'' \geq 0$ such that $c_j'+c_j''=O(d_{|j|}/|j|)$ as $|j| \to \infty$ and
\begin{equation}
\label{cond-arb-coeff}
\frac{1}{d_n}\sum_{j=0}^{n}c_j' \to a', \quad  \frac{1}{d_n}\sum_{j=0}^{n}c_j'' \to a'', \quad \frac{1}{d_n}\sum_{j=-n}^{0}c_j' \to -b', \quad \frac{1}{d_n}\sum_{j=-n}^{0}c_j'' \to -b'',
\end{equation}
for some $a',a'',b',b'' \in \bR$ with $|a'|+|b'|>0$ and $|a''|+|b''|>0$, where $d_n$ satisfies (\ref{rate-dn}) for some $1/\alpha \leq H<1$.
Let $a=a'-a''$ and $b=b'-b''$. Let $(a_n)_n$ be chosen to satisfy (\ref{def-an-alpha2}) or (\ref{def-an-alpha}) depending on whether $\alpha=2$ or $\alpha \in (0,2)$. Then
$$\frac{1}{d_n a_n} S_n(\cdot) \indistes \Lambda_{\alpha,H,a,b}(\cdot),$$
where $\indistes$ denotes the convergence in distribution specified by Definition \ref{def-star-conv}.
In particular,
$$\frac{1}{d_n a_n} S_n(\cdot) \To_{\cD}\Lambda_{\alpha,H,a,b}(\cdot) \quad \mbox{in} \quad (\bD[0,1],S).$$
\end{theorem}

\noindent {\bf Proof:} For any $i \in \bZ$, we define
$X_{i}'=\sum_{j \in \bZ}c_{j}'\xi_{i-j}$ and $X_{i}''=\sum_{j \in \bZ}c_{j}''\xi_{i-j}$.
Then $S_{n}(t)=S_{n}'(t)-S_{n}''(t)$, where $S_{n}'(t)=\sum_{i=1}^{[nt]}X_{i}'$ and
$S_{n}''(t)=\sum_{i=1}^{[nt]}X_{i}''$.

By Theorem \ref{th-pos-coeff},
$$\frac{1}{d_n a_n}S_{n}'(\cdot) \To_{\cD} \Lambda'(\cdot) \quad \mbox{in} \quad (\bD[0,1],M_1)$$
$$\frac{1}{d_n a_n}S_{n}''(\cdot) \To_{\cD} \Lambda''(\cdot) \quad \mbox{in} \quad (\bD[0,1],M_1),$$
where $\Lambda'(\cdot)=\Lambda_{\alpha,H,a',b'}(\cdot)$ and $\Lambda''(\cdot)=\Lambda_{\alpha,H,a'',b''}(\cdot)$.

By the linearity property of the stable integrals (page 117 of \cite{ST94}),
$$\Lambda'(t)-\Lambda''(t)=\int_{-\infty}^{\infty} [f_{\alpha,H,a',b'}(t,u)-f_{\alpha,H,a'',b''}(t,u)]dZ(u)=\Lambda_{\alpha,H,a,b}(t) \quad \mbox{a.s.}$$

The conclusion will follow from Theorem \ref{sum-conv}, once we prove that for any $0\leq t_1< \ldots <t_k \leq 1$,
\begin{equation}
\label{fin-dim-conv}
\frac{1}{d_n a_n}(S_{n}'(t_1), S_{n}''(t_1), \ldots,S_{n}'(t_k), S_{n}''(t_k) ) \To_{\cD} (\Lambda'(t_1), \Lambda''(t_1), \ldots, \Lambda'(t_k), \Lambda''(t_k)).
\end{equation}

Note that each of the components of the random vector on the left-hand side of (\ref{fin-dim-conv}) is a stochastic integral with respect to $Z_n$, in the sense of  (\ref{def-int-Zn}). More precisely,
$$\frac{1}{d_n a_n}S_{n}'(t)=\int_{-\infty}^{\infty}f_{n}'(t,u)dZ_n(u), \quad \frac{1}{d_n a_n}S_{n}''(t)=\int_{-\infty}^{\infty}f_{n}''(t,u)dZ_n(u)$$
with $$f_{n}'(t,u)=\frac{1}{d_n}\sum_{k=1-[nu]}^{[nt]-[nu]}c_k'
 \quad \mbox{and} \quad f_{n}''(t,u)=\frac{1}{d_n}\sum_{k=1-[nu]}^{[nt]-[nu]}c_k''.$$
The components of the random vector on the right-hand side of (\ref{fin-dim-conv}) are stochastic integrals with respect to $\{Z(t)\}_{t \in \bR}$. To prove (\ref{fin-dim-conv}), we apply Corollary 3.3 of \cite{KM88} to the functions:
$$f_n^1(u)=f_{n}'(t_1,u), f_n^2(u)=f_{n}''(t_1,u), \ldots,f_n^{2k-1}(u)=f_{n}'(t_k,u),$$ $$f_n^{2k}(u)=f_{n}''(t_k,u),$$
$$f^1(u)=f_{\alpha,H,a',b'}(t_1,u), f^2(u)=f_{\alpha,H,a'',b''}(t_1,u), \ldots,
f^{2k-1}(u)=f_{\alpha,H,a',b'}(t_k,u), $$ $$f^{2k}(u)=f_{\alpha,H,a'',b''}(t_k,u).$$

As in the proof of Theorem 5.1 of \cite{KM88}, one can show that for any $t \in [0,1]$ fixed, the functions $f_{n}'(t,\cdot)$ and $f_{n}''(t,\cdot)$ satisfy the following conditions:\\
(B1) $f_{n}'(t,u) \stackrel{{\rm c.c.}}{\To} f_{\alpha,H,a',b'}(t,u)$ and $f_{n}''(t,u) \stackrel{{\rm c.c.}}{\To} f_{\alpha,H,a'',b''}(t,u) $ a.e.($u$). Here, $\stackrel{{\rm c.c.}}{\To}$ denotes the continuous convergence, i.e. $f_n (u)\stackrel{{\rm c.c.}}{\To} f(u)$ at $u=u_0$ if $f_n(u_n) \to f(u_0)$ whenever $u_n \to u_0$.  \\
(B2) For any $T>0$ there exists $\beta>\alpha$ such that $$\sup_{n \geq 1} \int_{|u| \leq T} |f_{n}'(t,u)|^{\beta}\rho_n(du)<\infty \quad \mbox{and} \quad \sup_{n \geq 1} \int_{|u| \leq T} |f_{n}''(t,u)|^{\beta}\rho_n(du)<\infty,$$
where $\rho_n$ is the measure on $\bR$ defined by $\rho_n([0,u])=[nu]/n$.\\
(B3) There exists $\varepsilon>0$ such that
$$\lim_{T \to \infty}\limsup_{n \to \infty} \int_{|u|>T} (|f_{n}'(t,u)|^{\alpha-\varepsilon} +|f_{n}'(t,u)|^{\alpha+\varepsilon})d\rho_n(u)=0$$
$$\lim_{T \to \infty}\limsup_{n \to \infty} \int_{|u|>T} (|f_{n}''(t,u)|^{\alpha-\varepsilon} +|f_{n}''(t,u)|^{\alpha+\varepsilon}) d\rho_n(u)=0.$$

Therefore, the functions $f_n^{1}, \ldots, f_n^{2k}, f^1, \ldots, f^{2k}$ satisfy the conditions $(A1)'$, $(A2)'$ and $(A3)'$ of Corollary 3.3 of \cite{KM88}. Relation (\ref{fin-dim-conv}) follows. $\Box$

\vspace{3mm}

Our final result gives a set of conditions on the coefficients $(c_j)_{j \in \bZ}$, for which we can give an answer to the open problem mentioned in the introduction.

\begin{theorem}
\label{answer-open-problem}
Let $\alpha \in (0,2]$. Suppose that $c_j=c_j'-c_j''$ where $c_j',c_j'' \geq 0$,
$$\sum_{j \in \bZ}(c_j')^{\delta}<\infty \quad \mbox{and} \quad \sum_{j \in \bZ}(c_j'')^{\delta}<\infty$$
for some $0<\delta <\alpha$, $\delta \leq 1$. Let $A=\sum_{j \in \bZ}c_j$ and $(a_n)_n$ be chosen to satisfy (\ref{def-an-alpha2}) or (\ref{def-an-alpha}) depending on whether $\alpha=2$ or $\alpha \in (0,2)$. Then
$$\frac{1}{a_n}S_n(\cdot) \indistes AZ(\cdot) \quad \mbox{in} \quad (\bD[0,1],S),$$
where $\{Z(t)\}_{t \geq 0}$ is the process described in Remark \ref{def-Z-process} and $\indistes$ denotes the convergence in distribution specified by Definition \ref{def-star-conv}. In particular,
$$\frac{1}{a_n}S_n(\cdot) \To_{\cD} AZ(\cdot) \quad \mbox{in} \quad (\bD[0,1],S).$$
\end{theorem}

\noindent {\bf Proof:} The argument is similar to the one used in the proof of Theorem \ref{main-result}. The idea is to first prove the result for non-negative coefficients $(c_j)_{j \in \bZ}$ (similarly to Theorem \ref{th-pos-coeff}), and then apply Theorem \ref{sum-conv}. To prove the finite dimensional convergence of the pair $(S_n'/a_n,S_n''/a_n)$, we use again Corollary 3.3 of \cite{KM88} (as in the proof of Theorem 5.3 of \cite{KM88}). $\Box$

\section{Examples and Simulations}
\label{section-examples}

In this section we consider several examples to illustrate the results presented in Section \ref{section-result}.
We let $\xi,(\xi_i)_{i \in \bZ}$ be i.i.d. random variables with a symmetric Pareto distribution with parameter $\alpha \in (0,2]$, i.e. $\xi$ has density
$$f(x)=\frac{1}{2}\alpha x^{-\alpha-1}1_{(1,\infty)}(x)+\frac{1}{2}\alpha (-x)^{-\alpha-1}1_{(-\infty,-1)}(x).$$
If $\alpha>1$, $E(\xi)=0$. Note that $|\xi|$ has a Pareto distribution with parameter $\alpha$, i.e. $|\xi|$ has density $g(x)=\alpha x^{-\alpha-1} 1_{(1,\infty)}(x)$.

When $\alpha<2$, we choose $a_n$ such that $nP(|\xi|>a_n) \to 1$. A suitable choice is
$$a_n=\inf\{x>0;P(|\xi| \leq x) \geq 1-n^{-1}\}=G^{-1}(1-n^{-1})=n^{1/\alpha}$$
where $G(x)=P(|\xi| \leq x)=1-x^{-\alpha}$ and $G^{-1}(y)=(1-y)^{-1/\alpha},y \in [0,1]$. In this case,
$\{Z(t)\}_{t \geq 0}$ is an $\alpha$-stable L\'evy motion with $Z(1) \sim S_{\alpha}(\sigma,0,0)$ and the scale parameter $\sigma$ is given by (\ref{def-sigma}) with $C=1$.

When $\alpha=2$, we choose $a_n$ as the largest root of the equation $x^2=2n \ln x$, since the function $U(x)=E(|\xi|^2 1_{\{|\xi| \leq x\}})=2 \ln x$ is slowly varying, and hence $\xi$ is in the domain of attraction of the normal law. (However, note that $E(\xi^2)=\infty$.) In this case, $\{Z(t)\}_{t \geq 0}$ is a standard Brownian motion. More precisely, we used the formula $a_n=\exp\{-\frac{1}{2}W_{-1}(-1/n)\}$ where $W_{-1}$ is the second branch of the Lambert $W$ function (see e.g. \cite{CGHJK96}).

For the simulations, we used $$\xi_i=F^{-1}(U_i)=-(2U_i)^{-1/\alpha}1_{[0,1/2]}(U_i)+[2(1-U_i)]^{-1/\alpha}1_{(1/2,1]}(U_i)$$ where $(U_i)_{i \in \bZ}$ are i.i.d. random variables with a uniform distribution on $[0,1]$, and $F^{-1}$ is the generalized inverse of the c.d.f. $F$ of $\xi$:
$$F(x)=\left\{
\begin{array}{ll} 2^{-1}(-x)^{-\alpha} & \mbox{if $x \leq -1$} \\
2^{-1} & \mbox{if $x \in (-1,1]$} \\
1-2^{-1}x^{-\alpha} & \mbox{if $x \geq 1$}
\end{array} \right.$$
We used the truncated series $X_i^{N}=\sum_{|j| \leq N}c_{j}\xi_{i-j}$ as an approximation for $X_i$, and the corresponding partial sum sequence $S_n^N(\cdot)=\sum_{i=1}^{[n \cdot]}X_i^N$ as an approximation for $S_n(\cdot)$, if $N$ is large. We considered $n=1000$ and $N=50$. Therefore, $a_n=95.4883$ when $\alpha=2$.

\begin{remark}
\label{conv-plots}
{\rm (i) Usually, to illustrate a functional limit theorem by the convergence of plots, one uses the automatic scaling done by the computer, as explained in Section 1.2.4 of \cite{whitt02}. For instance, for illustrating Donsker's theorem for the partial sum sequence $S_n(\cdot)=\sum_{i=1}^{[n\cdot]}\xi_i$ associated with i.i.d. random variables $(\xi_i)_{i \geq 1}$ with mean 0 and variance 1, one would like to plot the step function determined by the points $(k/n,S_k/\sqrt{n})$ for $k=1, \ldots,n$, where $S_n=S_n(1)$. To do this and guarantee that the plot fits in the available space, the computer automatically shifts the values $T_{k,n}=S_k/\sqrt{n}$ by $\min_{k \leq n}T_{k,n}$ and then scales them by a factor equal to the range $r=\max_{k \leq n}T_{k,n}-\min_{k \leq n}T_{k,n}$.
Since the map $plot:\bD[0,1] \to \bD[0,1]$ defined by $plot(x)=(x-\inf_{t \in [0,1]} x(t))/(\sup_{t \in [0,1]} x(t)-\inf_{t \in [0,1]} x(t))$ is $J_1$-continuous, by the continuous mapping theorem,
\begin{equation}
\label{conv-plots-Donsker}
plot\left(\frac{S_n(\cdot)}{\sqrt{n}} \right) \To_{\cD} plot(Z(\cdot)) \quad \mbox{in} \ (\bD[0,1],J_1),
\end{equation}
where $\{Z(t)\}_{t \in [0,1]}$ is a standard Brownian motion.
Therefore, the convergence of the plots is preserved, despite the automatic shift-and-scale operation.

(ii) In our case, we had to avoid this automatic shift-and-scale operation since the map $plot$ is not $S$-continuous. However, the map $plot_{a,b}:\bD[0,1] \to \bD[0,1]$ defined by $plot_{a,b}(x)=(x-a)/(b-a)$ is $S$-continuous for any $a,b \in \bR$ with $a<b$. Therefore, for the simulations in Example \ref{ex-Th-main-result}  below we chose ourselves the range for the values $T_{k,n}^N=S_k^N/(d_n a_n)$ with $k=1, \ldots,n$, where $S_n^N=S_n^N(1)$ (for $n$ and $N$ fixed). More precisely, we
forced a shift-and-scale of these values in a more ``realistic'' range $[x_{\rm min},x_{\rm max}]$. This range was obtained by simulating $(T_{k,n}^N)_{k \leq n}$ a large number $M$ of times, recording the values $T_{\rm min}=\min_{k \leq n}T_{k,n}^N$ and $T_{\rm max}=\max_{k \leq n}T_{k,n}^N$ each time, and taking $x_{\rm min}$ (respectively $x_{\rm max}$) as the 10\%-quantile of the $M$ values $T_{\rm min}$ (respectively the
90\%-quantile of the $M$ values $T_{\rm max}$). We considered $M=75$. The same procedure was employed in Example \ref{ex-Th-answer-open-problem} for $T_{k,n}^N=S_{k}^N/a_n$, using the 15\%-quantile, respectively the 85\%-quantile.
}
\end{remark}

\begin{remark}
{\rm (i) The procedure that we explained in Remark \ref{conv-plots}.(i) provides some reasonable justification for the illustration of Donsker's theorem using the plots. However, there is a small problem with this justification, since relation (\ref{conv-plots-Donsker}) gives a convergence in distribution, and not an a.s. convergence. In reality, what we illustrate with this procedure is a convergence of plots which holds {\em almost every time} we perform the simulation. This turns out to be a consequence of the Koml\'os-Major-Tusn\'ady strong approximation result (Theorem 2 of \cite{KMT2}):
$$\sup_{t \in [0,1]} \left|\frac{S_n(t)}{\sqrt{n}}-Z(t)\right|=o(n^{-\ve}) \quad \mbox{a.s.}$$
with $\ve=2^{-1}-p^{-1}$, assuming that $E|\xi_i|^{p}<\infty$ for some $p>3$.

(ii) In our case, the fact that the plots resemble the trajectories of the desired limit process suggests that it might be possible to prove some strong approximation results which would parallel the results presented here. We do not investigate this problem here.
}
\end{remark}


\begin{example}[Illustration of Theorem \ref{th-pos-coeff}]
\label{ex-Th-pos-coeff}
{\rm We consider the coefficients:
$$c_j=j^{-\gamma} \ \mbox{for all} \ j \geq 1 \quad \mbox{and} \quad c_j=0 \ \mbox{for all} \ j \leq 0.$$
We have two cases:
(i) $H>1/\alpha$; (ii) $H=1/\alpha$.

To illustrate (i), we let $\alpha^{-1}<\gamma<1$ and $d_n=n^{1-\gamma}=n^{H-1/\alpha}$ where
\begin{equation}
\label{def-H}
H=\frac{1}{\alpha}+1-\gamma.
\end{equation}
The fact that $\gamma>\alpha^{-1}$ guarantees that (\ref{cond-coeff}) holds. Since
$$\sum_{j=1}^{n}c_j=\sum_{j=1}^{n}j^{-\gamma} \sim \int_1^{n+1} x^{-\gamma}dx \sim \frac{1}{1-\gamma}n^{1-\gamma},$$
the hypotheses of Theorem \ref{th-pos-coeff} are verified with $L_2(n)=1$, $a=(1-\gamma)^{-1}$ and $b=0$.
Figure \ref{figure-th3-7i} gives an approximation of the sample path of the process $\{\Lambda_{\alpha,H,a,0}(t)\}_{t \in [0,1]}$ obtained for $\gamma=0.75$; hence $a=4$. The picture on the left was obtained using $\alpha=1.5$ (the limit process is a LFSM with $H=0.92$), while the picture on the right corresponds to $\alpha=2$ (the limit process is the FBM of index $H=0.75$, multiplied by $a/C_H=4.28$, where $C_H$ is given by (\ref{definition-CH}).) 
\begin{figure}
\begin{center}
\includegraphics[scale=0.25]{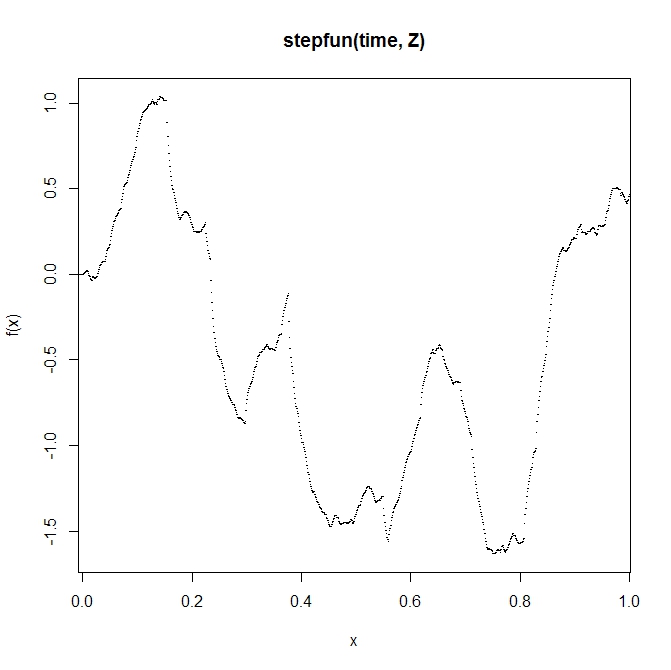}
\includegraphics[scale=0.25]{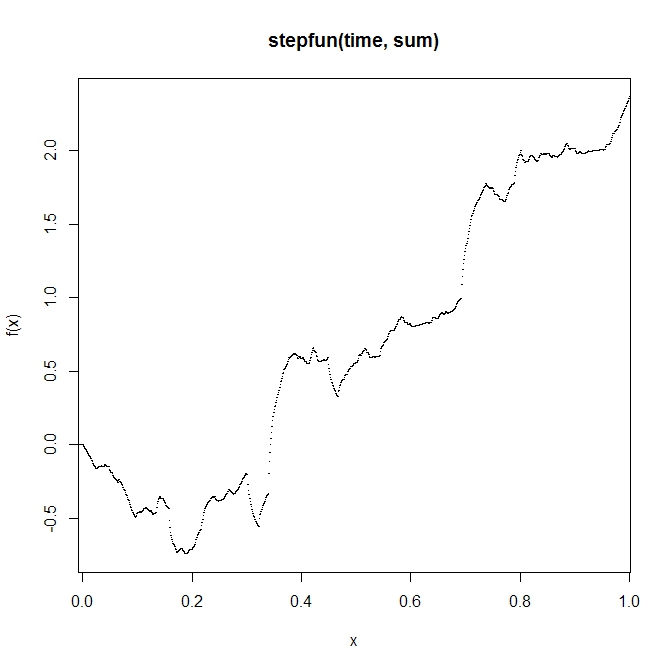}
\end{center}
\caption{Illustration of Theorem 3.7 for $H>\alpha^{-1}$: $\alpha=1.5$ (left), $\alpha=2$ (right)}
\label{figure-th3-7i}
\end{figure}

To illustrate (ii), we let $\gamma>1$ (hence $\gamma>\alpha^{-1}$). We consider $d_n=1$. The hypotheses of Theorem \ref{th-pos-coeff} are verified since
$$\frac{1}{d_n}\sum_{j=1}^{n}c_j \to \sum_{j \geq 1}j^{-\gamma}=\zeta(\gamma)=:a$$
where $\zeta$ is the Riemann-zeta function.
Figure \ref{figure-th3-7ii} gives an approximation of the sample path of the process $\{aZ(t)\}_{t \in [0,1]}$ obtained for $\gamma=4$, and hence $a=\zeta(4)=\pi^4/90=1.08$ (see p.807 of \cite{AS65}). The picture on the left was obtained using $\alpha=1.5$ (the limit is an $\alpha$-stable L\'evy motion), while for the picture on the right we used $\alpha=2$ (the limit is the Brownian motion multiplied by $a$).
\begin{figure}
\begin{center}
\includegraphics[scale=0.25]{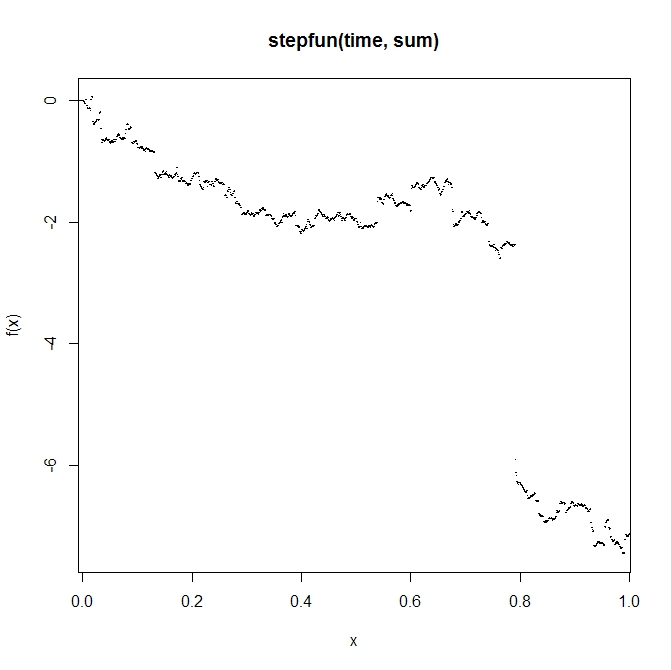}
\includegraphics[scale=0.25]{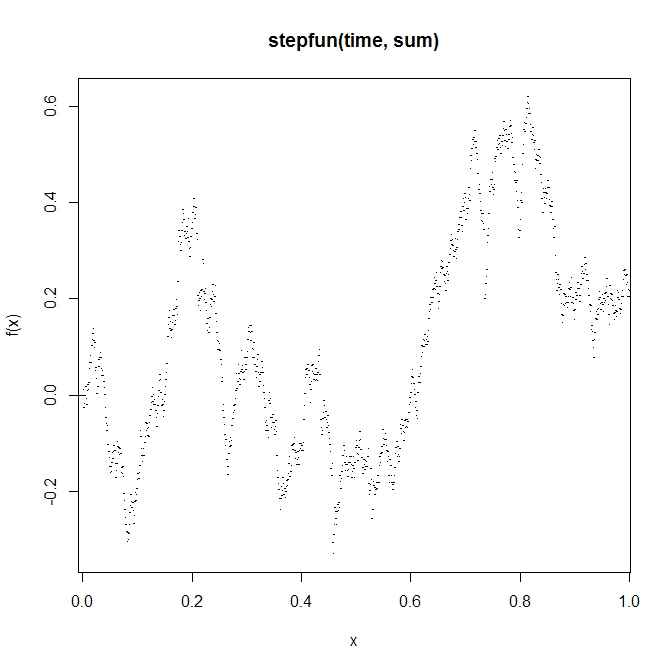}
\end{center}
\caption{Illustration of Theorem 3.7 for $H=\alpha^{-1}$: $\alpha=1.5$ (left), $\alpha=2$ (right)}
\label{figure-th3-7ii}
\end{figure}

In both figures, the shift-and-scale of the values $S_k^N/(d_n a_n)$ was done automatically by the computer, since the convergence holds in the $M_1$ topology.
}
\end{example}

\begin{example}[Illustration of Theorem \ref{main-result}]
\label{ex-Th-main-result}
{\rm We consider the coefficients:
\begin{equation}
\label{sign-coeff}
c_j=\left\{
\begin{array}{ll}
k_1 j^{-\gamma} & \mbox{if $j$ is even} \\
-k_2 j^{-\gamma} & \mbox{if $j$ is odd}
\end{array} \right. \mbox{for all} \ j \geq 1 \quad \mbox{and} \quad c_j=0 \ \mbox{for all} \ j \leq 0,$$
for some $k_1,k_2>0$. We write $c_j=c_j'-c_j''$ where
$$c_j'=\left\{
\begin{array}{ll}
k_1 j^{-\gamma} & \mbox{if $j$ is even} \\
0 & \mbox{if $j$ is odd}
\end{array} \right. \quad \mbox{and} \quad
c_j''=\left\{
\begin{array}{ll}
0 & \mbox{if $j$ is even} \\
k_2 j^{-\gamma} & \mbox{if $j$ is odd}
\end{array} \right. \quad \mbox{for} \ j \geq 1
\end{equation}
and $c_j'=c_j''=0$ for $j \leq 0$.
We have two cases: (i) $H>1/\alpha$; (ii) $H=1/\alpha$.

To illustrate (i), we let $\alpha^{-1}<\gamma<1$ and $d_n=n^{1-\gamma}=n^{H-1/\alpha}$ where $H$ is given by (\ref{def-H}). Since $\gamma>\alpha^{-1}$, (\ref{cond-coeff}) holds. Note that
\begin{eqnarray*}
\sum_{j=1,j \ {\rm even}}^{n}j^{-\gamma} &=&\sum_{k=1}^{[n/2]}(2k)^{-\gamma}\sim 2^{-\gamma} \frac{1}{1-\gamma}\left[ \frac{n}{2}\right]^{1-\gamma} \sim \frac{1}{2(1-\gamma)}n^{1-\gamma} \\
\sum_{j=1,j \ {\rm odd}}^{n}j^{-\gamma} &=&\sum_{j=1}^{n}j^{-\gamma} -\sum_{j=1,j \ {\rm even}}^{n}j^{-\gamma}\sim \frac{1}{2(1-\gamma)}n^{1-\gamma}.
\end{eqnarray*}
Therefore,
$$\frac{1}{d_n}\sum_{j=1}^{n}c_j' \to \frac{k_1}{2(1-\gamma)}=:a' \quad \mbox{and} \quad
\frac{1}{d_n}\sum_{j=1}^{n}c_j'' \to \frac{k_2}{2(1-\gamma)}=:a''.$$ Hence, $a=a'-a''=(k_1-k_2)/[2(1-\gamma)]$. Note that $a=0$ if $k_1=k_2$.
\begin{figure}
\begin{center}
\includegraphics[scale=0.25]{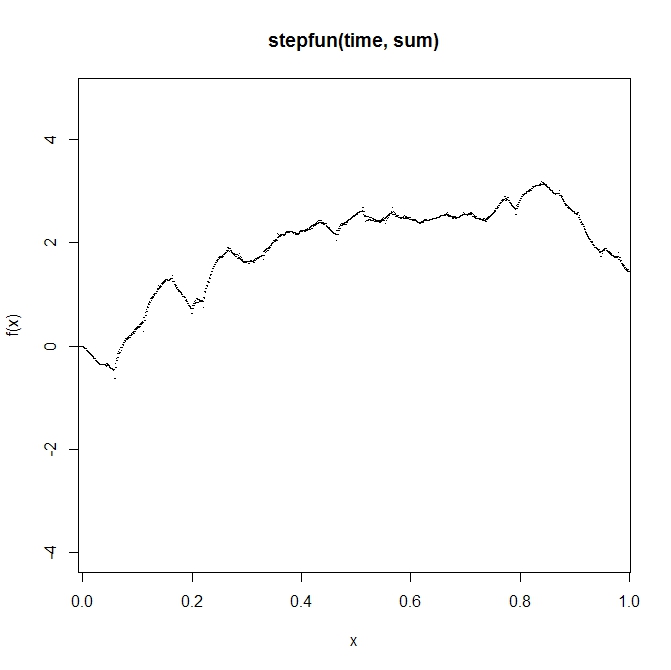}
\includegraphics[scale=0.25]{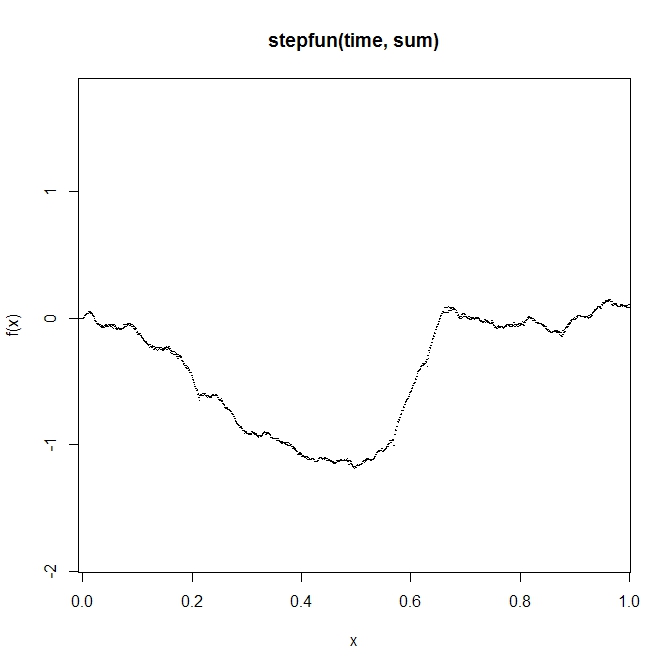}
\end{center}
\caption{Illustration of Theorem 3.9 for $k_1=3,k_2=1,H>\alpha^{-1}$: $\alpha=1.5$ (left), $\alpha=2$ (right)}
\label{figure-th3-9i-1}
\end{figure}

Figure \ref{figure-th3-9i-1} gives an approximation for a sample path of the process $\{\Lambda_{\alpha,H,a,0}(t)\}_{t \in [0,1]}$ obtained for $k_1=3,k_2=1, \gamma=0.75$ (hence $a=4$) and $\alpha=1.5$ (left), respectively $\alpha=2$ (right). On the left, the limit process is a LFSM with $H=0.92$; on the right, the limit is the FBM of index $H=0.75$ multiplied by $a/C_H=4.28$.
The plot of the values $S_{k}^N/(d_n a_n)$ was performed in the range $[-4.03,4.83]$ (left), respectively $[-1.86,1.75]$ (right).

The same values $\gamma=0.75$ and $\alpha=1.5$ were used for Figure \ref{figure-th3-9i-2}
but with $k_1=k_2=1$. In this case, {\em the limit is the zero process} since $a=0$. In the the picture on the left, we imposed the shift-and-scale operation in the interval $[-0.62,0.52]$, while in the picture on the right we used the automatic shift-and-scale performed by the computer. Therefore, the picture on the right is a blow-up of the picture on the left. (Note the small values on the $y$-axis in the picture on the right.)
\begin{figure}
\begin{center}
\includegraphics[scale=0.25]{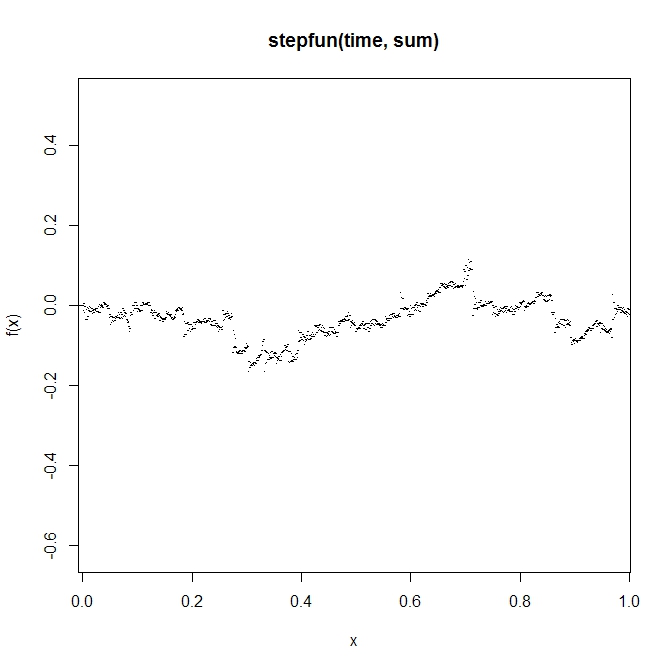}
\includegraphics[scale=0.25]{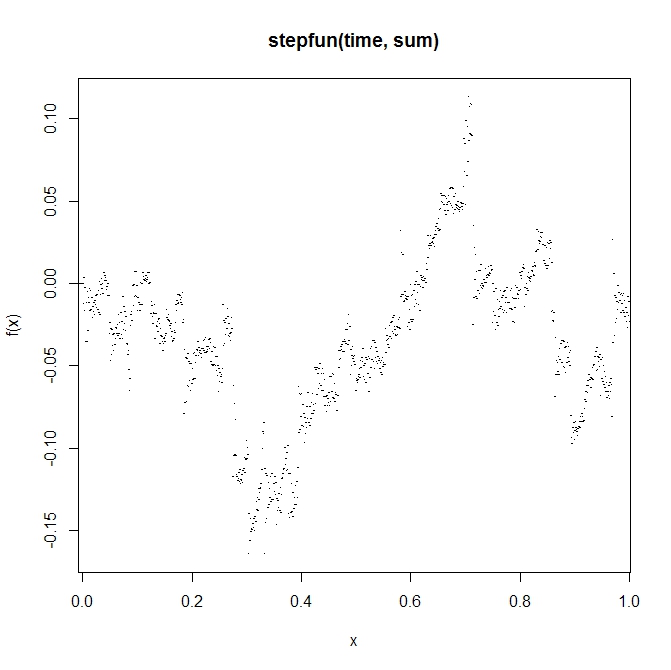}
\end{center}
\caption{Illustration of Theorem 3.9 for $k_1=k_2=1$ ($a=0$) and $H>\alpha^{-1}$: shift-and-scale in the range [-0.62,0.52] (left), automatic shift-and-scale (right)}
\label{figure-th3-9i-2}
\end{figure}

To illustrate (ii), we let $\gamma>1$, $d_n=1$ and $H=1/\alpha$. Then
$$\frac{1}{d_n}\sum_{j=1}^{n}c_j' \to k_1 \sum_{j \geq 1, j \ {\rm even}}j^{-\gamma}=:a' \quad \mbox{and} \quad \frac{1}{d_n}\sum_{j=1}^{n}c_j'' \to k_2 \sum_{j \geq 1, j \ {\rm odd}}j^{-\gamma}=:a''.$$
In this case, $a=a'-a''=\sum_{j \geq 1}c_j$, but $a \not=0$ when $k_1=k_2$.

Figure \ref{figure-th3-9ii-1} gives an approximation for a sample path of the process $\{a Z(t)\}_{t \in [0,1]}$ obtained for $k_1=k_2=1,\gamma=4$ and $\alpha=1.5$ (left), respectively $\alpha=2$ (right). In this case, $a=\sum_{j \geq 1}(-1)^{j}j^{-4}=-7\pi^4/720=-0.95$ (see p.807 of \cite{AS65}). The picture on the left is an approximation of an $\alpha$-stable L\'evy motion, while the picture on the right is an approximation of the Brownian motion multiplied by $a$. The plot of the values $S_{k}^N/(d_n a_n)$ was performed in the range $[-3.49,4.16]$ (left), respectively $[-1.55,1.45]$ (right).
\begin{figure}
\begin{center}
\includegraphics[scale=0.25]{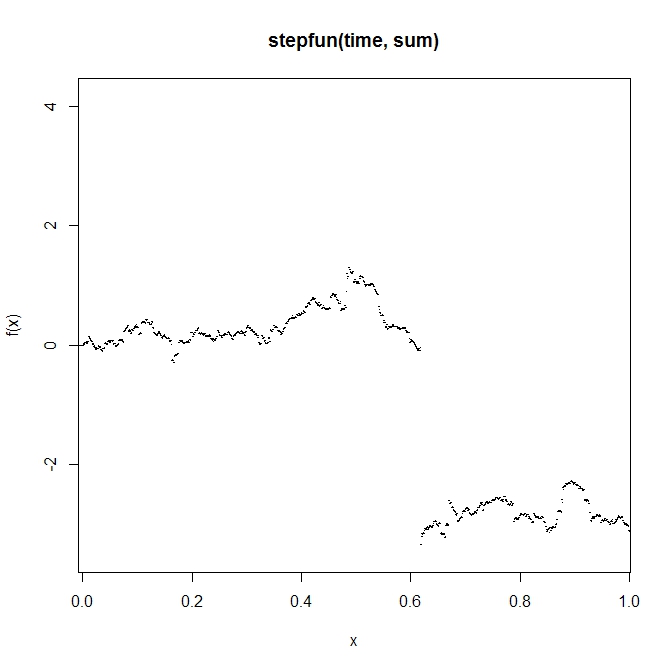}
\includegraphics[scale=0.25]{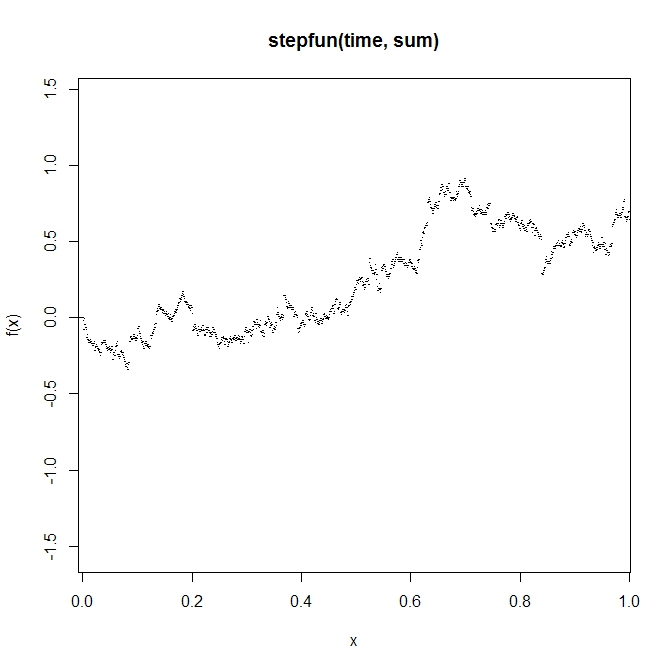}
\end{center}
\caption{Illustration of Theorem 3.9 for $k_1=k_2=1,H=\alpha^{-1}$: $\alpha=1.5$ (left), $\alpha=2$ (right)}
\label{figure-th3-9ii-1}
\end{figure}

}
\end{example}

\begin{example}[Illustration of Theorem \ref{answer-open-problem}]
\label{ex-Th-answer-open-problem}
{\rm We assume that the coefficients $(c_j)_{j \in \bZ}$
are given by (\ref{sign-coeff}). In order that $\sum_{j \in \bZ}|c_j|^{\delta}<\infty$ for some $0<\delta<\alpha,\delta \leq 1$, we need $\gamma>\max(\alpha^{-1},1)$. We let
$$A=\sum_{j \in \bZ}c_j=k_1 \sum_{j \geq 1, j \ {\rm even}}j^{-\gamma}-k_2\sum_{j \geq 1, j \ {\rm odd}}j^{-\gamma} $$

To illustrate the result, we consider $k_1=k_2=1$ and $\gamma=4$ (hence $A=-0.95$).
Figure \ref{figure-th3-10}
gives and approximation for a sample path of the process $\{AZ(t)\}_{t \in [0,1]}$ in the case $\alpha=0.8$ (left), respectively $\alpha=1$ (right). The plot of the values $S_{k}^N/a_n$ was performed in the range $[-4.38,2.38]$ (left), respectively $[-2.32,4.23]$ (right).
\begin{figure}
\begin{center}
\includegraphics[scale=0.25]{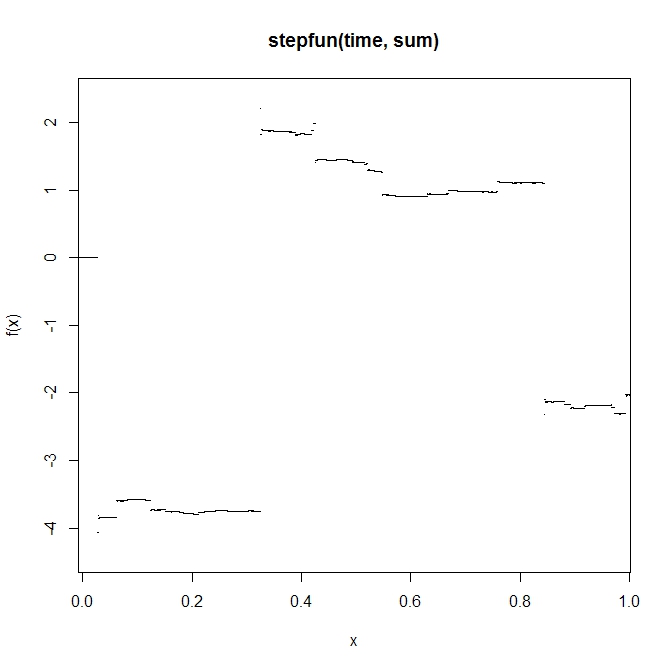}
\includegraphics[scale=0.25]{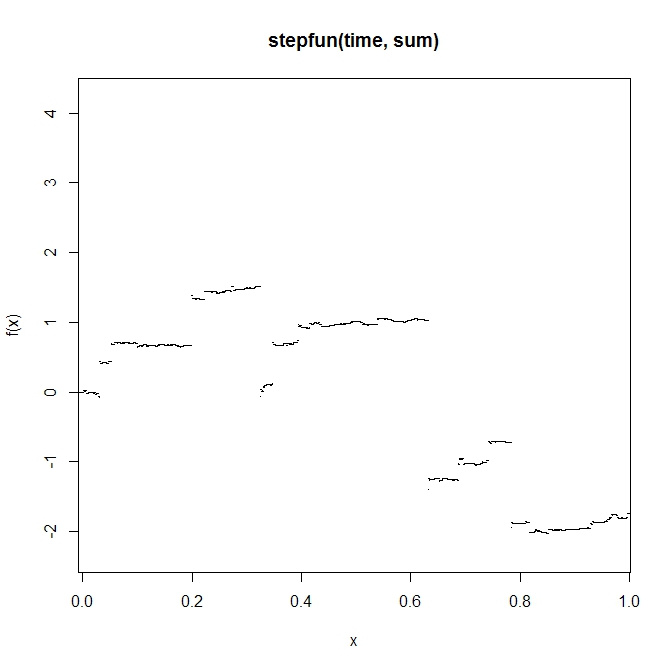}
\end{center}
\caption{Illustration of Theorem 3.10: $\alpha=0.8$ (left), $\alpha=1$ (right)}
\label{figure-th3-10}
\end{figure}

}
\end{example}

\appendix
\section{Some auxiliary results}

\begin{lemma}
\label{lem-m1-bound}
Let $x \in \bD[0,1]$ be arbitrary. For any $0 \leq s \leq u <v \leq t \leq 1$,
$$|x(u)-x(v)| \leq 2|x(s)-x(t)|+H(x(s),x(u),x(t))+H(x(s),x(v),x(t)).$$
\end{lemma}

\noindent {\bf Proof:} We consider only the case $x(s) \leq x(t)$, the case $x(t) <x(s)$ being similar. We claim that for any $s \leq u \leq t$,
\begin{equation}
\label{lemA-step1}
|x(u)-x(s)| \leq |x(t)-x(s)|+H\big(x(s),x(u),x(t) \big).
\end{equation}
To see this, we consider three cases.
If $x(u)<x(s)$, then $|x(u)-x(s)|=H\big(x(s),x(u),x(t)\big)$. If $x(s) \leq x(u) \leq x(t)$ then $|x(u)-x(s)| \leq |x(t)-x(s)|$. Finally, if $x(u)>x(t)$ then $|x(u)-x(s)| \leq H\big(x(s),x(u),x(t)\big)+|x(t)-x(s)|$.

The conclusion follows using relation (\ref{lemA-step1}) for $u$ and $v$, and the fact that $|x(u)-x(v)| \leq |x(u)-x(s)|+|x(v)-x(s)|$. $\Box$

\begin{lemma}
\label{lem-m1-number}
Let $x \in \bD[0,1]$ be arbitrary. For $0 \leq s<t \leq 1$, define
$$\beta=\sup_{s \leq u<v<w \leq t}H(x(u),x(v),x(w)).$$
If $\eta>2\beta$ then
$$N_{\eta}(x;[s,t]) \leq \frac{2|x(t)-x(s)|+\beta}{\eta-\beta},$$
where $N_{\eta}(x;[s,t])$ denotes the number of $\eta$-oscillations of $x$ in the interval $[s,t]$.
\end{lemma}

\noindent {\bf Proof:} Let $s \leq t_1<t_2 \leq t_3<t_4 \leq \ldots \leq t_{2N-1}<t_{2N} \leq t$ be such that
$$|x(t_{2k})-x(t_{2k-1})|>\eta \quad \mbox{for all} \ k=1,\ldots,N.$$
Assume first that $x(t_2)-x(t_1)>\eta$. We claim that:
$$x(t_3)\geq x(t_2)-\beta \quad \mbox{and} \quad  \quad x(t_4)-x(t_3)>\eta.$$
To see this, suppose that $x(t_3)<x(t_2)-\beta$. Then the distance between $x(t_2)$ and the interval with endpoints $x(t_1)$ and $x(t_3)$ is greater than $\beta$, which is a contradiction. Hence $x(t_3)\geq x(t_2)-\beta$. On the other hand, if we assume that $x(t_4)-x(t_3)<-\eta$, we obtain that
$$x(t_1)=x(t_1)-x(t_2)+x(t_2)-x(t_3)+x(t_3)<-\eta+\beta+x(t_3)<x(t_3)-\beta,$$
which means that the distance between $x(t_3)$ and the interval with endpoints $x(t_1)$ and $x(t_4)$ is greater than $\beta$, again a contradiction.

Repeating this argument, we infer that:
$$x(t_{2k})-x(t_{2k-1})>\eta, \quad \mbox{for all} \quad k=1, \ldots,N$$
and
$$x(t_{2k+1})-x(t_{2k})>-\beta \quad \mbox{for all} \quad k=1, \ldots, N-1.$$
Taking the sum of these inequalities, we conclude that:
\begin{equation}
\label{lemB-step1}
x(t_{2N})-x(t_1)>N\eta-(N-1)\beta=N(\eta-\beta)+\beta.
\end{equation}

On the other hand, by Lemma \ref{lem-m1-bound}, we have:
\begin{equation}
\label{lemB-step2}
|x(t_{2N})-x(t_1)|\leq 2|x(t)-x(s)|+2\beta.
\end{equation}

Combining (\ref{lemB-step1}) and (\ref{lemB-step2}), we obtain that
$$N \leq \frac{2|x(t)-x(s)|+\beta}{\eta-\beta},$$
which is the desired upper bound.
$\Box$

\end{document}